\documentclass[11pt]{article}

\usepackage{babel} 

\newenvironment{poliabstract}[1]
  {\renewcommand{\abstractname}{#1}\begin{abstract}}
  {\end{abstract}}

\usepackage[top= 2cm,bottom=2cm,left=2.00cm,right=2.00cm]{geometry}

\usepackage{amsfonts,amsmath,floatflt,graphicx,multicol}
\usepackage{hyperref}

\usepackage[super,sort&compress]{natbib}
\bibpunct{}{}{,}{s}{}{,}
\renewcommand{\d}{{\rm d}}
\newcommand{\e}{{\rm e}}
\renewcommand{\i}{{\rm i}}
\newcommand{\D}{{\rm D}}

\newcommand{\PD}[2]{\frac{\partial #1}{\partial #2}}
\newcommand{\FD}[2]{\frac{\d #1}{\d #2}}

\renewcommand{\vec}[1]{\mathbf{#1}}

\DeclareMathSymbol{\ZSet}{\mathalpha}{AMSb}{"5A}
\DeclareMathSymbol{\RSet}{\mathalpha}{AMSb}{"52}
\DeclareMathSymbol{\CSet}{\mathalpha}{AMSb}{"43}

\title{\vspace*{-1.5cm}Revisiting the Haken Lighthouse Model}

\author{S. Coombes\footnote{School of Mathematical Sciences, University of Nottingham, Nottingham.
NG7 2RD, UK. \newline {\tt email:} stephen.coombes@nottingham.ac.uk}}
\date{}

\begin{document}

\maketitle
\begin{poliabstract}{Abstract} 
Simple spiking neural network models, such as those built from interacting integrate-and-fire (IF) units, exhibit rich emergent behaviours but remain notoriously difficult to analyse, particularly in terms of their pattern-forming properties. In contrast, rate-based models and coupled phase oscillators offer greater mathematical tractability but fail to capture the full dynamical complexity of spiking networks. To bridge these modelling paradigms, Hermann Haken --- the pioneer of Synergetics --- introduced the Lighthouse model, a framework that provides insights into synchronisation, travelling waves, and pattern formation in neural systems.

In this work, we revisit the Lighthouse model and develop new mathematical results that deepen our understanding of self-organisation in spiking neural networks. Specifically, we derive the linear stability conditions for phase-locked spiking states in Lighthouse networks structured on graphs with realistic synaptic interactions ($\alpha$-function synapses) and axonal conduction delays. Extending the analysis on graphs to a spatially continuous (non-local) setting, we develop a variant of Turing instability analysis to explore emergent spiking patterns. Finally, we show how localised spiking bump solutions --- which are difficult to mathematically analyse in IF networks --- are far more tractable in the Lighthouse model and analyse their linear stability to wandering states.

These results reaffirm the Lighthouse model as a valuable tool for studying structured neural interactions and self-organisation, further advancing the synergetic perspective on spiking neural dynamics.
\end{poliabstract}

\begin{poliabstract}{Preface} 
I was fortunate to meet Hermann Haken at the \textit{STOCHAOS} (Stochastic and Chaotic Dynamics in the Lakes) conference in Ambleside, UK in 1999.  At that time I was a post-doc working with Paul Bressloff at Loughborough, UK and we were investigating the dynamics of spiking integrate-and-fire (IF) networks.  At STOCHAOS, I gave one of my first presentations on the work with Paul, and Hermann was very kind to come and talk to me later at the meeting and give encouragement for further work.  Prior to that I only knew his work on Synergetics from afar, so was very pleased to have this direct interaction.  At this meeting I heard for the first time about his Lighthouse model \cite{Haken2000b}, though did not pursue it any further as I did not have the foresight to see the advantage over studying IF networks.  A few years later in 2004 I visited Carson Chow at NIH in Bethesda, and we talked about searching for simpler models than IF, as by then both of us had realised it wasn't always such an easy model to work with.  Here, Carson suggested a spike-based model that after a while I realised was  non other than the Lighthouse model.  Despite out slight disappointment that we weren't the first to ponder such a model, we did manage to extoll its virtues by showing how one could explicitly construct spiking bump network solutions (though we did not treat their stability) \cite{Chow2006}.  I am delighted to revisit the Lighthouse model in this paper and obtain a suite of new results for waves, bumps, and patterns, that I hope Hermann would have found interesting.
\end{poliabstract}

\newpage

\section{Introduction}

A central goal of theoretical neuroscience is to construct tractable models of spiking neuronal networks. Achieving this requires a single-neuron model capable of generating discrete spikes of activity --- action potentials --- that, when embedded in a synaptic network, can give rise to the rich and complex behaviours observed in biological neural systems. Despite considerable efforts, a comprehensive understanding of network dynamics has remained elusive for both biophysically detailed conductance-based neuron models and simplified integrate-and-fire (IF) models. This challenge largely stems from the lack of a robust mathematical framework to describe the collective dynamics of spiking networks.
To date, most theoretical progress has been made using firing rate models, as reviewed in \cite{Coombes2014}, which average out spiking activity and thus fail to capture key features such as spike-train correlations.  Other approaches, including those based on spike density formulations \cite{Gerstner2014}, often rely on IF-type neurons with global coupling and pulsatile interactions, incorporating a predefined stochastic firing rate. While these methods offer some analytical traction, they remain limited in scope and do not fully reflect the complexity of real spiking network dynamics.

The Lighthouse model for action potential coupled neuronal networks, originally introduced by Haken at the end of the 20th century, offers a promising bridge between spiking and firing rate descriptions, and mixes in ideas from phase-oscillator networks for good measure \cite{Haken2000a,Haken2000b}.  In the regime of slow synaptic interactions, the Lighthouse model reduces to familiar neural mass models of the type exemplified by that of Wilson \& Cowan \cite{Wilson1972}. In contrast, under fast synaptic dynamics, it exhibits a range of complex behaviours characteristic of spiking networks.  Extensive work by Haken in his book ``Brain Dynamics: Synchronization and Activity Patterns in Pulse-coupled Neural Nets with Delays and Noise"  \cite{Haken2002}, has characterised its phase-locked solutions and emphasised the model's tractability at the network level.  Nonetheless, to date it has received very little further attention in the theoretical neuroscience community.   Here, we aim to redress this imbalance by emphasising its analytical tractability in scenarios that range from discrete to continuous neural networks, covering synchrony, waves, and patterns (both localised and globally spatially periodic), emphasising solution structure and stability in terms of spike times.

In Sec.~\ref{Sec:Lighthouse} we introduce the Haken Lighthouse model of a spiking neuronal network, posed on a graph, with structured synaptic interactions triggered by action potential arrival.  Next in Sec.~\ref{Sec:SyncandStab} we show that, when the weight connection matrix satisfies a row sum constraint, synchrony is a naturally emergent solution.  Moreover, a (map based) analysis of the network firing times is shown to be suitable for deriving the conditions for linear stability.  Here, we also develop a so-called (flow based) \textit{Master Stability Function} approach that allows one to determine the linear stability of the synchronous solution for arbitrary networks in a computationally pragmatic fashion.  The extension of the Lighthouse model to spatially continuous spiking systems is treated in Sec.~\ref{Sec:Continuum}.  Here, we also highlight that, in the limit of slow synaptic dynamics, the Lighthouse model reduces to a Wilson-Cowan-like rate-based model, thereby providing a formal connection between spiking and rate-based approaches.
With the inclusion of space-dependent axonal delays in the spiking model we show that it readily supports periodic travelling waves, and determine their speed (as a function of period) in Sec.~\ref{Sec:TW}, and stability in Sec.~\ref{Sec:Turing}.  Here, we also show how to determine the Turing instability of the synchronous state to globally periodic spatio-temporal patterns.  The treatment of spatially localised \textit{bumps} is treated next, in Sec.~\ref{Sec:Bump}, in the case that the nonlinearity of the Lighthouse model is taken to be a Heaviside.  We show that bumps can be built from propagating spiking waves and that our analysis (construction and stability) generalises that proposed by Amari for firing rate neural fields.  At various points throughout the paper, we also consider the limit of slow synaptic processing where firing rate models can be recovered.  This is useful for emphasising the more novel behaviours of the full spiking model, as well as showing that theoretical results for spiking solutions recover those of rate based reductions (in the limit of slow synapses).  Finally, in Sec.~\ref{Sec:Discussion}, we recapitulate the main points of this paper and discuss potential areas for further exploration.

\section{The Lighthouse model\label{Sec:Lighthouse}}

\begin{figure}
\centering
\includegraphics[width=0.6\textwidth]{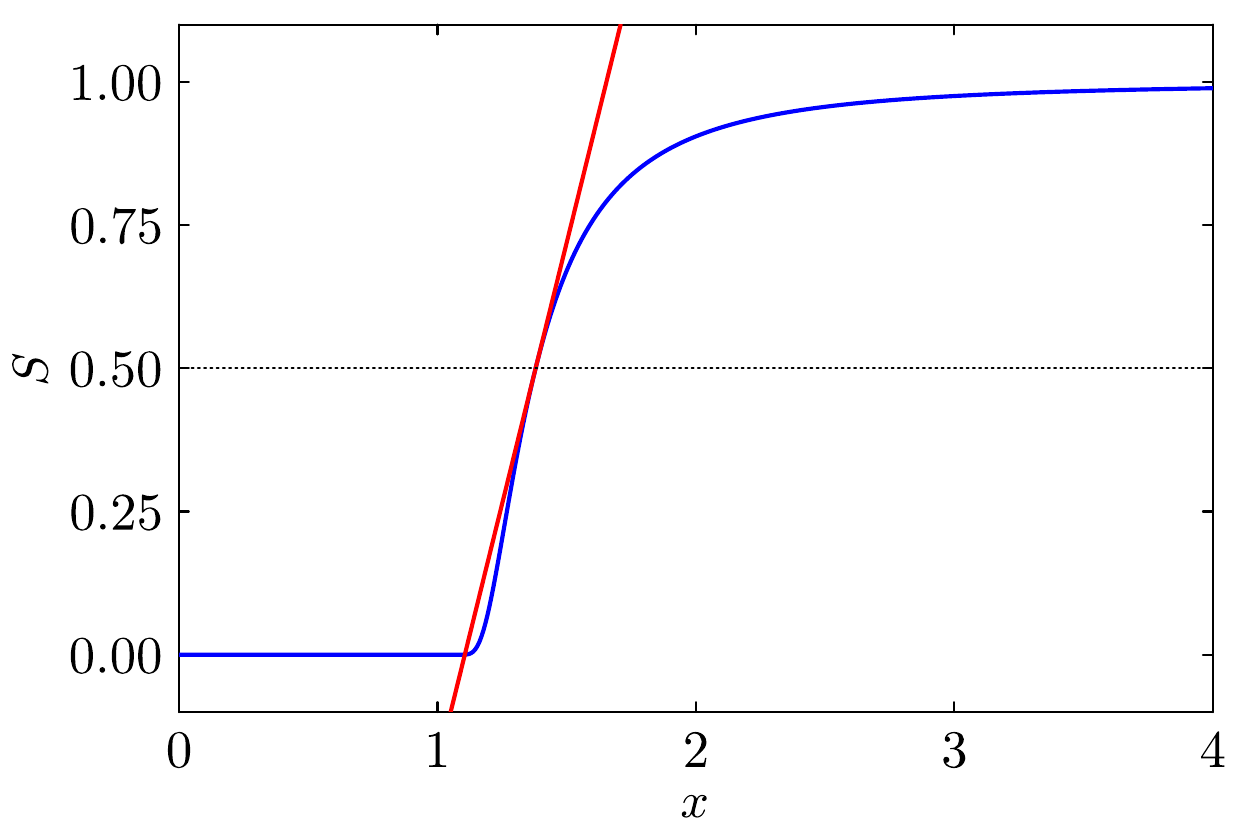}
\caption{
\label{Fig:SPlot}
A plot of the nonlinear function $S$ (in blue) given by equation (\ref{S}) and its mid-range linear approximation $S_L$ (in red).  
Parameters: $r=0.1$ and $h=1$.
}
\end{figure}

Let us introduce a number of nodes indexed by $i=1,2,\ldots, N$.  At each node $i$ we attach scalar dynamical variables $\theta_i(t)$ and $\psi_i(t)$, representing the state of the system and the input delivered to node $i$ at time $t >0$ respectively.  The Haken Lighthouse model posed on a graph with connection strengths $w_{ij} \in \RSet$ between nodes $i$ and $j$ can then be written as a dynamical system in the form
\begin{equation}
\FD{}{t} \theta_i (t)= S \left ( \psi_i(t) \right ) , \qquad \psi_i(t)  = \sum_{j=1}^N w_{ij} s_j(t - \tau_{ij}), \qquad s_i(t) = \sum_{m \in \ZSet} \eta(t-T_i^m) .
\label{LHmodel}
\end{equation}
Here, $\tau_{ij} \geq 0$ is a communication delay between nodes $i$ and $j$, $S$ is a non-negative function, and $\eta(t)$ encapsulates the shape of a post-synaptic response to the arrival of an action potential.  The $m$th time of action potential generation at node $i$ is denoted $T_i^m$.  
The state variable $\theta_i$ may be regarded as the lift of an angle and the firing times are defined by the condition
\begin{equation}
\lim_{\delta \rightarrow 0_+} \theta_i(T_i^m - \delta) \! \! \! \!\mod 2\pi  = 2 \pi .
\end{equation} 
Haken required the function $S$ to have the following properties: $S(x)$ is equal to zero for $x$ smaller than a threshold $h$ and then it increases in a quasi-linear fashion until it saturates.
In \cite{Haken2002} Haken chose to work with the Naka--Rushton formula (Hill formula with a cut-off)
though here we prefer a slightly different representation of $S$ and one convenient for analysis, namely the $C^\infty$ function:
\begin{equation}
S(x) = \exp \left ( - \frac{r}{(x-h)^2} \right ) H(x - h), \qquad r>0, 
\label{S}
\end{equation}
where $H$ is a Heaviside function.  A graph of $S(x)$ and a linear approximation in its mid-range is shown in Fig.~\ref{Fig:SPlot}.  Note that $ \lim_{r \rightarrow 0} S(x) = H(x-h)$.
It is pertinent to point out the Haken originally formulated two version so the Lighthouse model.  In one version, whenever the input is below threshold or drops below threshold the phase immediately resets to zero.  In the second version the phase is not reset.  We shall only consider the second version here.

For the nonlinear choice of $S$ given by (\ref{S}) the interpretation of $\theta_i$ as the lift of an angle means that if $\psi_i(t) > h$ for all $t$ then the Lighthouse model can be thought of as a phase oscillator with a non-uniform rotation.  This rotation becomes uniform in the limit $r \rightarrow 0$, much like the behaviour of an actual lighthouse emitting beams of light from a system of lamps and lenses that serves as a beacon for navigational aid for ships at sea.  From a more neural perspective, the Lighthouse model accommodates a natural refractory process, in the sense that it takes at least a time $2 \pi$ to reach the firing condition.  However, in contrast to many other phase-oscillator networks for which the uncoupled system reduces to a network of independent oscillators, the Lighthouse model will only behave in this way if $h<0$.   In this case, an isolated Lighthouse neuron does not have any novel response to external inputs in the sense that its (infinitesimal) phase response curve is flat.

Haken originally treated an exponential synapse of the form $\eta(t) = \alpha \e^{-\alpha t} H(t)$.  Given that this has a jump discontinuity at the origin which may not be biophysically realistic we prefer to work with the continuous $\alpha$-function given by
\begin{equation}
\eta(t) = \alpha^2 t \exp(-\alpha t) H (t) .
\label{alpha}
\end{equation}
Note that the theory we develop below is relatively agnostic to the actual form of $\eta$ as long as it is causal and has a well defined Fourier transform.  For the rest of the paper we shall assume that this is the case with the further restriction that $\eta$ is normalised such that $\int_{0}^\infty \d t \, \eta(t) = 1$.

\section{Synchrony and its stability\label{Sec:SyncandStab}}

In \cite{Haken2000a,Haken2000b,Haken2002} Haken studied the synchrony of a pair of interacting Lighthouse oscillators with $S(x)$ a linear function and in the absence of delays.  Here we consider larger networks with a single fixed delay $\tau_{ij} = \tau$ for all $(i,j)$ and focus on networks with a row sum constraint $\sum_{j} w_{ij} = \Gamma$ for all $i$.  This latter condition guarantees the existence of a synchronous solution.  

\subsection{Emergent period}

\begin{figure}[htbp]
\centering
\includegraphics[width=0.6\textwidth]{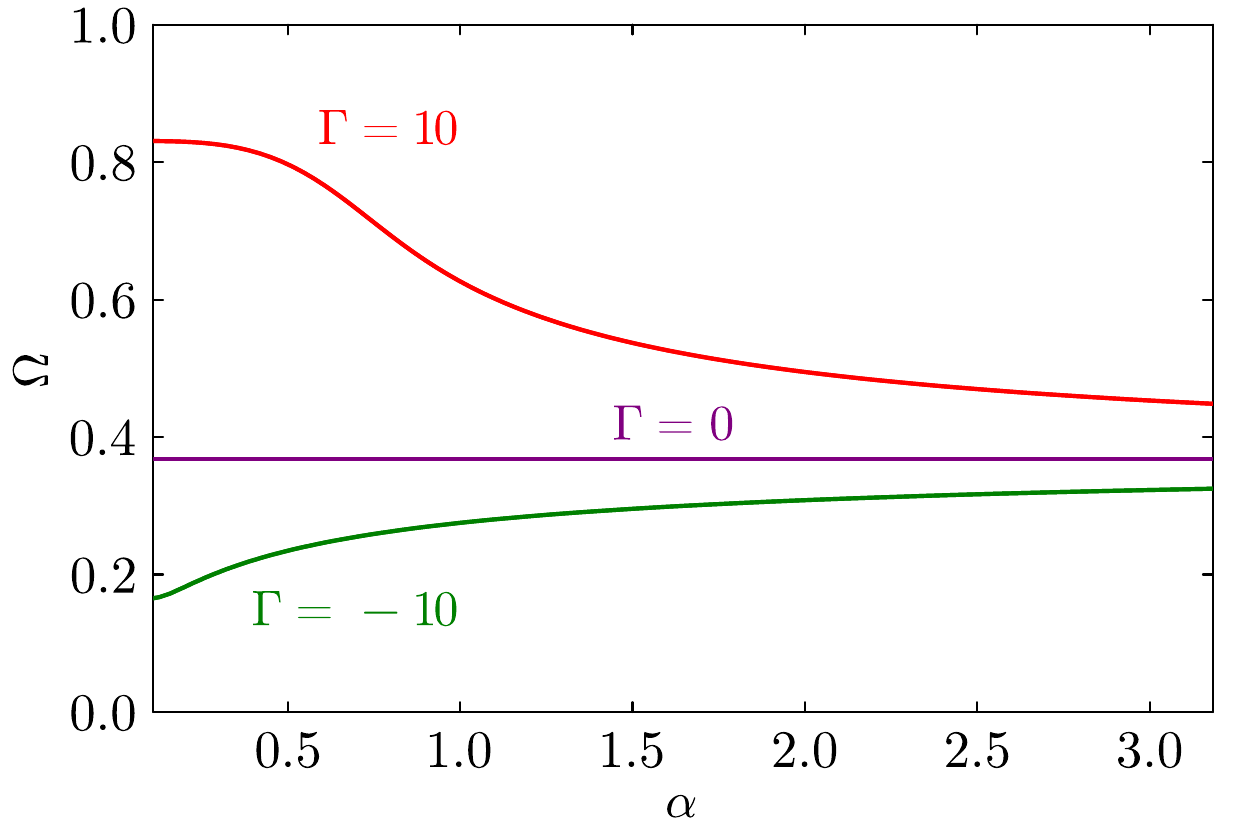}
\caption{The emergent frequency $\Omega = 2 \pi/T$ of the (nonlinear) Lighthouse network model with a network structure such that $\sum_{j} w_{ij} = \Gamma$ for all $j$ (row sum constraint) as a function of the synaptic rate parameter $\alpha$.  Parameters: $r=1$, $h=-1$, and $\tau=0$.
\label{Fig:Period}
}
\end{figure}

We define a synchronous network solution by $\theta_i(t) = \theta(t)$ for all $i$ with firing times $T_i^m = mT$ for some period $T$.  In this case $\psi_i(t) = \psi(t)$ for all $i$ where
$\psi(t) = \Gamma P(t)$, where $P(t)$ is a $T$-periodic function given by
\begin{equation}
P(t) = \sum_{m \in \ZSet} \eta (t-\tau - mT) .
\end{equation}
This may be equivalently represented by a Fourier series of the form $P(t) = \sum_n P_n \e^{\i \omega_n t }$ with $\omega_n = 2 \pi n/T$
\begin{equation}
P_n = \frac{1}{T} \widehat{\eta}(\omega_n) \e^{ -\i \omega_n \tau}, \quad \widehat{\eta}(\omega) = \int_{-\infty}^ \infty \d t \, \eta(t) \e^{-\i \omega t} .
\end{equation}
Here $\widehat{\eta}$ is recognised as the Fourier transform of $\eta$.  
For the $\alpha$-function (\ref{alpha}) we have that
\begin{equation}
\widehat{\eta} (\omega) = \frac{1}{(1+\i \omega/\alpha)^2} .
\end{equation}
We may now integrate (\ref{LHmodel}) to obtain
\begin{equation}
\theta(t) -\theta(0) = \int_0^t S \left ( \Gamma P(s) \right ) \d s .
\end{equation}
The emergent period of oscillation is then determined from the constraint $\theta(T) - \theta(0)  = 2 \pi$.  
For the case of a \textit{balanced} network  with $\Gamma = 0$ we have the simple result that
\begin{equation}
T = \frac{2 \pi}{S(0)}.
\end{equation}
We note that balance occurs naturally in networks with a \text{graph Laplacian structure}:  $w_{ij} = - a_{ij} + \delta_{ij} \sum_k a_{ik}$, for which it is always the case that $\sum_j w_{ij} = 0$ for all $i$.  
For the more general case when $\Gamma \neq 0$ the period can be found numerically solving the implicit equation $F(T)=0$, where $F(T) = 2 \pi  -\int_0^T \d s \,  S \left ( \Gamma P(s) \right )$.
However, if the Lighthouse model operates in a linear regime then more explicit progress can be made.  Consider for example the mid-range restriction $S \rightarrow S_L$, where $S_L(x) = 1/2 +S'(x_{1/2}) (x-x_{1/2})$ and $x_{1/2}$ is defined by $S(x_{1/2}) = 1/2$.  Equivalently we may write this in the form $S_L(x) = \gamma x - \Theta$, where $\gamma = S'(x_{1/2})>0$ and $\Theta = S'(x_{1/2}) x_{1/2} - 1/2 \in \RSet$.  In this case the emergent period satisfies the equation
\begin{equation}
2 \pi  =  \gamma \Gamma  \int_0^T P(s) \d s - \Theta T  .
\label{}
\end{equation}
Since $P(t)$ is a $T$-periodic function $T^{-1} \int_0^T P(s) \d s = P_0$ and using the fact that $P_0=T^{-1}$ we arrive at a formula for the period as
\begin{equation}
T = \frac{\gamma \Gamma - 2 \pi}{\Theta} ,
\label{Tlinear}
\end{equation}
which essentially recovers a result first obtained by Haken for two oscillators \cite{Haken2000a}[eqn (125)].  
we note that for $\Theta >0$ we have the natural condition $\gamma \Gamma > 2 \pi$ (and vice versa).
For a balanced network we have the simple result that $T = 2 \pi /|\Theta|$ for $\Theta <0$
In Fig.~\ref{Fig:Period} we show a plot of the emergent frequency $\Omega = 2 \pi/T$ for the full nonlinear function (\ref{S}) obtained by the numerical solution of $2 \pi - \int_0^T S \left ( \Gamma P(s) \right ) \d s = 0$ as a function of $\alpha$ for various choices of $\Gamma$.

\subsection{Linear stability\label{sec:linstab}}

Here, as in Haken \cite{Haken2000a,Haken2000b,Haken2002}, we shall consider the case that $S$ is linear and fix $S=S_L$.  To progress we borrow some ideas from \cite{Bressloff2000} (for IF networks) and introduce a perturbed trajectory $\widetilde{\theta}_i(t)$ that results from perturbations in the firing times.  We denote the perturbed firing times by $\widetilde{T}_i^m$.  We further introduce the state and time deviations $\delta\theta_i(t) = \widetilde{\theta}_i(t) - \theta(t)$ and $\delta T_i^m = \widetilde{T}_i^m - mT$ respectively.  We may relate these state and time deviations by using the firing condition:
\begin{equation}
 \theta (m T)  = \widetilde{\theta}_i(\widetilde{T}_i^m) .
 \label{fire-fire}
\end{equation}
A Taylor expansion of (\ref{fire-fire}) for small deviations gives the result
\begin{equation}
\delta T_i^m = - \frac{\delta \theta_i (mT) }{\dot{\theta} (T)} .
\label{dT}
\end{equation}
The evolution equation for the perturbed trajectory is given by 
\begin{equation}
\FD{}{t} \widetilde{\theta}_i (t) = \gamma  \sum_{j} w_{ij} \sum_{m \in \ZSet} \eta (t -\tau - \widetilde{T}_j^m) - \Theta \nonumber \\
\end{equation}
Integrating this from $\widetilde{T}_i^m$ to $\widetilde{T}_i^{m+1}$ gives
\begin{align}
\widetilde{\theta}_i(\widetilde{T}_i^{m+1}) &- \widetilde{\theta}_i(\widetilde{T}_i^m) =
\gamma \sum_{j} w_{ij} \int_{\widetilde{T}_i^m}^{\widetilde{T}_i^{m+1}} \d \, t \sum_{p \in \ZSet} \eta (t -\tau - \widetilde{T}_j^p) - \Theta \left (\widetilde{T}_i^{m+1} - \widetilde{T}_i^m \right ) \nonumber \\
& \simeq (\gamma \Gamma - \Theta T ) +(\gamma \Gamma P(T) - \Theta) (\delta T_i^{m+1} - \delta T_i^m) + \gamma \sum_j w_{ij} \sum_p G_p [ \delta T_i^m - \delta T_j^{m-p}] ,
\end{align}
where 
\begin{equation}
G_p = \int_0^T \d t \, \eta' (t-\tau+pT) .
\end{equation}
Using the fact that $\dot{\theta}(t) = \gamma \Gamma P(t) - \Theta$ and remembering (\ref{dT}) gives
\begin{equation}
\delta \theta_i((m+1)T) - \delta \theta_i(mT) = \frac{\gamma}{\dot{\theta}(T)} \sum_j w_{ij} \sum_p G_p [\delta \theta_j((m-p)T) - \delta \theta_i(mT) ] .
\end{equation}
This linear difference equation has a discrete spectrum that can be found by taking solutions of the form $\delta \theta_i(mT) = u_i \e^{m \lambda}$ where $\lambda \in \CSet$.
This yields the  eigen-system:
\begin{equation}
\dot{\theta}(T) \left ( \e^{\lambda} - 1 \right ) u_i = \gamma \sum_j w_{ij} \left[  G(\lambda) u_j - G(0) u_i \right ] , \qquad G(\lambda) = \sum_p G_p \e^{-p \lambda} .
\end{equation}

A computationally useful form for $G(\lambda)$ can be obtained by using a Fourier integral representation of $\eta'$ as:
\begin{equation}
\eta'(t) = \frac{1}{2 \pi} \int_{-\infty}^\infty \d k \, (\i k) \widehat{\eta} (k) \e^{\i k t} .
\end{equation}
Hence,
\begin{align}
G(\lambda) 
&= \sum_{p \in \ZSet}  \frac{1}{2 \pi} \int_{-\infty}^\infty \d k \,  \left ( \e^{\i k T} - 1 \right ) \widehat{\eta} (k) \e^{\i k (-\tau+pT)}   \e^{- p \lambda} \nonumber \\
& = \frac{1}{T} \sum_{n \in \ZSet} \widehat{\eta} (\omega_n -\i \lambda/T) \e^{-\i (\omega_n -\i \lambda/T) \tau} (\e^{\lambda}-1) , \label{Gl}
\end{align}
where we have made use of the Dirac comb identity $2 \pi \sum_{n \in \ZSet} \delta(x-2 \pi n) = \sum_{p \in \ZSet} \e^{\i p x}$.  We note from (\ref{Gl}) that $G(0) = 0$, and this result could have been anticipated earlier by noting that 
\begin{equation}
G(0) = \sum_{p \in \ZSet} G_p = \int_0^T \d t \, \sum_{p \in \ZSet} \eta' (t-\tau+pT)  = \sum_{p \in \ZSet} \left [ \eta (T-\tau+pT) - \eta (0-\tau+pT) \right ] = P(T) - P(0) = 0 . 
\end{equation}

If $w$ is diagonalisable with eigenvalues $\widehat{w}_\mu$, $\mu = 0, \ldots, N-1$ then the characteristic equations for the spectrum can be written $ ( \e^{\lambda} - 1 )\mathcal{E}_\mu(\lambda) = 0$, where
\begin{equation}
\mathcal{E}_\mu(\lambda) =   \dot{\theta}(T ) - \gamma \widehat{w}_\mu \mathcal{G}(\lambda)  ,
\label{E}
\end{equation}
and
\begin{equation}
\mathcal{G}(\lambda) = \frac{1}{T} \sum_{n \in \ZSet} \widehat{\eta} (\omega_n -\i \lambda/T) \e^{-\i (\omega_n -\i \lambda/T) \tau} .
\end{equation}

We note that $\lambda=0$ is always an eigenvalue (and arises because of time-translation invariance of the synchronous oscillatory solution). 
The synchronous network solution will be stable provided that all other eigenvalues (defined by zeros of $\mathcal{E}_\mu(\lambda) =0$) have real part negative.  A natural way to numerically obtain the eigenvalues from $\mathcal{E}_\mu (\lambda)$ is to decompose $\lambda = \nu + \i \omega$ and plot the zero levels sets of $\text{Re} \, \mathcal{E}_\mu (\nu + \i \omega) = 0$ and $\text{Im} \, \mathcal{E}_\mu (\nu + \i \omega) = 0$.  The eigenvalues can then be found by the intersection of the two level set curves.  The synchronous solution is stable if all eigenvalues reside in the left hand complex plane.

For $\lambda \in \RSet$ then $\mathcal{G} (\lambda) \in \RSet$ and for a matrix $w$ with only real eigenvalues it may be that a single real eigenvalue can cross the origin in the complex plane.  This defines a static bifurcation leading to the emergence of non-synchronous behaviour, which could be in the form of a phase-locked network state or a network state with heterogeneous frequencies.  The conditions for such a real instability can be found by solving $\dot{\theta}(T) = \gamma \widehat{w}_\mu \mathcal{G}(0)$ for some system parameter.  
For $\lambda \in \CSet$ we may seek a dynamic instability by considering $\lambda = \i \omega$ for $\omega \in \RSet  \backslash 0$.  In this case, decomposing $\widehat{w}_\mu = \widehat{w}_R + \i \widehat{w}_I$ and introducing $\mathcal{G}_R(\omega) = \text{Re} \, \mathcal{G}(\i \omega)$ and $\mathcal{G}_I (\omega)= \text{Im} \, \mathcal{G}(\i \omega)$ we may equate the real and imaginary parts of $\mathcal{E}_\mu (\i \omega)= 0$ to obtain the pair of equations
\begin{align}
\label{Dynamic1}
\dot{\theta} (T) \cos (\omega) &= 1 + \gamma  \left [\widehat{w}_R \mathcal{G}_R(\omega) - \widehat{w}_I \mathcal{G}_I(\omega) \right ] \\
\label{Dynamic2}
\dot{\theta} (T) \sin (\omega) &= \gamma  \left [\widehat{w}_R \mathcal{G}_I(\omega) + \widehat{w}_I \mathcal{G}_R(\omega) \right ] .
\end{align}
Solving this pair of equations for $(\omega, a) = (\omega_c, a_c)$, where $a$ is any parameter in the model (such as $\gamma$, $\Theta$, $\Gamma$, $\alpha$ or $\tau$) fixes a pair of complex conjugate eigenvalues $\pm \i \omega_c$ on the the imaginary axis.  If all other eigenvalue for the choice $a = a_c$ have negative real part and the pair on the imaginary axis move into the right hand complex plane under variation of $a$ then a Hopf instability can arise.   This would lead to the formation of periodic solutions, that is, ones with a time-dependent firing rate.

\begin{figure}
\centering
\includegraphics[width=0.6\textwidth]{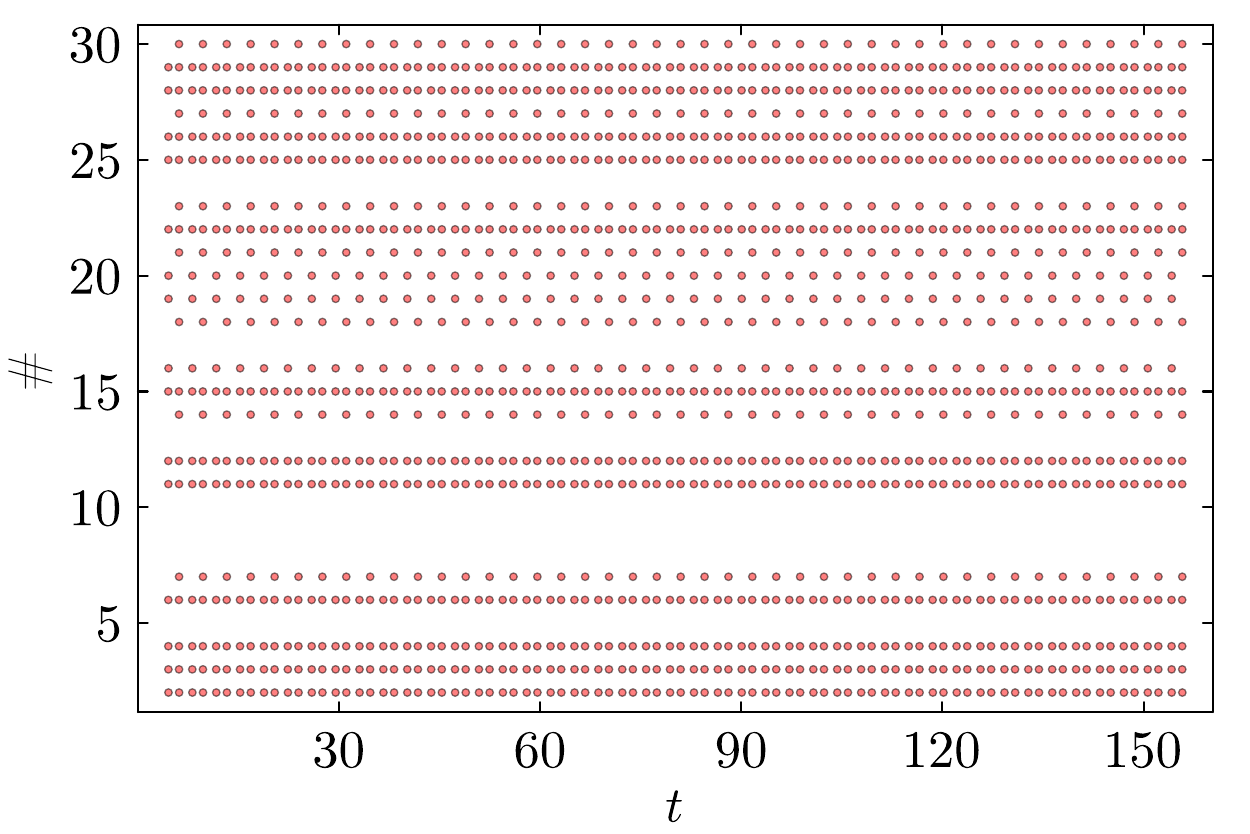}
\caption{A raster plot of spike times illustrating \textit{oscillator death} (arising from the instability of a synchronous state) in a globally coupled network with $w_{ij} =  (\Gamma+1) \delta_{ij} - N^{-1}$.  Open circles denote times of firing.
Parameters: $N=30$, $\Gamma=1$, $\gamma=\pi$, $\Theta=-1$, $\alpha = 5$, and $\tau=0$.
\label{Fig:OscillatorDeath}
}
\end{figure}

By way of example, consider a globally coupled network with $w_{ij} =  (\Gamma+1) \delta_{ij} - N^{-1}$.
$\Gamma$ is an eigenvalue of $w$ with eigenvector $(1,1,\ldots,1)$, and $\Gamma + 1$ is an $N-1$ degenerate eigenvalue with eigenvector $(u_1,u_2,\ldots,u_N)$ where $\sum_j u_j= 0$.  Beyond an instability one would expect to see the excitation of an inhomogeneous state (that can be predicted from the eigenvector pattern with $\sum_j u_j= 0$).  For $\gamma >0$ one can look for a critical value of $\alpha$ such that  a dynamic instability (as described above) is reached.  By solving equations (\ref{Dynamic1}) and (\ref{Dynamic2}) we may determine a pair $(\omega, \gamma) = (\omega_c, \gamma_c)$ for fixed $\alpha$.  Under variation in $\alpha$ we may then generate
a critical curve $\gamma = \gamma_c(\alpha)$.  For $\gamma > \gamma_c(\alpha)$ the system may be unstable to \textit{oscillator death}, whereby one or more of the oscillators ceases to fire.  An example of such behaviour, arising from the instability of a synchronous state, is shown in Fig.~\ref{Fig:OscillatorDeath} using a raster plot of spike times.

\subsubsection{Slow synapse reduction\label{sec:slow}}

For slow synapses, in the sense that $\alpha T/(2 \pi) \ll 1$, we need only keep the $n=0$ terms in the Fourier representations used above (since the Fourier coefficients all fall of rapidly).
In this case $P(t) = P_0 = 1/T$, so that $T$ satisfies the implicit equation $T = 2 \pi/S(\Gamma/T)$ for the full nonlinear Lighthouse model.  For the linearised model the period is still given by (\ref{Tlinear}) (as it is independent of $\alpha$).  Noting that $\dot{\theta}(T)$ reduces to $\dot{\theta}(T)= \gamma \Gamma/T - \Theta$ we may use (\ref{Tlinear}) to obtain the relation
\begin{equation}
\dot{\theta}(T) = \frac{2 \pi}{T} .
\label{slow_period}
\end{equation}
This gives the result that the rotational frequency of the oscillator at a firing event is equivalent to the firing frequency of the synchronous population.

For the linear stability of the linear Lighthouse model the characteristic equation simplifies to:
\begin{equation}
\frac{\e^{\lambda \tau/T}}{\widehat{\eta}(-\i \lambda /T)} = \frac{\gamma \widehat{w}_\mu}{2 \pi} ,
\label{lambdaslow}
\end{equation}
or equivalently
\begin{equation}
\left ( 1 + \frac{\lambda}{\alpha T} \right )^2 =\frac{\gamma \widehat{w}_\mu}{2 \pi}\e^{-\lambda \tau/T}.
\label{lambdaslowa}
\end{equation}
Consider for example the case that $w$ has real eigenvalues such that $\max_{\mu} \widehat{w}_\mu >0$  and $\min_{\mu} \widehat{w}_\mu < 0$.  A static bifurcation will occur at a critical value of $\gamma_c$ for which $1 = \gamma_c  \max_{\mu} \widehat{w}_\mu/(2 \pi )$. To determine the conditions for a delay induced dynamic instability we set $\lambda = \i \omega$ in (\ref{lambdaslowa}) and equate real and imaginary parts to obtain the pair of equations
\begin{equation}
1 -\left ( \frac{\omega}{\alpha T} \right )^2  = 2 \pi \gamma \widehat{w} \cos (\omega \tau/T) , \qquad 2  \left ( \frac{\omega}{\alpha T} \right )  = -2 \pi \gamma  \widehat{w}  \sin(\omega \tau/T) .
\end{equation}
Using the fact that $\cos^2 \theta + \sin^2 \theta = 1$ we obtain $\omega = \alpha T \sqrt{\pm \gamma \widehat{w}/(2 \pi) - 1}$.  Thus a necessary condition for a dynamic instability is 
$-\gamma  \min_{\mu} \widehat{w}_\mu/(2 \pi) \geq 1$ whilst preserving the condition $\gamma  \max_{\mu} \widehat{w}_\mu/(2 \pi ) <1$ (to avoid a real instability).
From the definition of $\tan \theta$ we obtain the formula $\tan(\omega \tau/T) = 2 \omega/(\alpha T)/(1-(\omega/(\alpha T))^2)$, for which the the smallest non-zero positive solution occurs when $0< \omega \tau/T < \pi/2$.  We observe that as $\tau$ increase from zero the largest possible solution for $\omega$ is $\pi T/(2 \tau)$.  Given that $-\gamma_c \min_{\mu} \widehat{w}_\mu /(2\pi) = (\omega/(\alpha T))^2+1$, we may obtain the approximate condition for a dynamic instability as
\begin{equation}
-\gamma  \min_{\mu} \widehat{w}_\mu/(2 \pi) - 1 = \left ( \frac{\pi T}{2 \tau} \right )^2 ,
\end{equation}
which would predict an emergent frequency $\omega \approx (\alpha T^2/(4 \tau)) 2 \pi$, and thus a period that scales as $4 \tau$.

%
%

\subsection{An equivalent MSF formulation}
\label{sec:MSF}

For the case of zero delays ($\tau=0$) the Master Stability Function (MSF) is a tool from dynamical network science that can assess the stability of the synchronous oscillatory state \cite{Pecora1998} in an arbitrary network with graph Laplacian coupling and smooth dynamics.  It can be readily extended to certain nonsmooth systems with balance, as recently reviewed in \cite[Ch.~7]{Coombes2023}.  We can apply this general framework to the linearised Lighthouse model with the introduction of the state vector $x_i = (\theta_i,s_i,u_i)$ and by writing the network dynamics in the equivalent form
\begin{equation}
\FD{x_i}{t} = A {x}_i + \gamma \sum_{j=1}^N w_{ij} {F}({x}_j) +b , \qquad i=1, \ldots, N,
\end{equation}
where $F(x_i) = (s_i,0,0)$, and
\begin{equation}
A= \begin{bmatrix}
0 & 0 & 0 \\
0 & -\alpha & \alpha \\
0 & 0 & -\alpha
\end{bmatrix},
\quad
b= -\begin{bmatrix}
\Theta\\
0\\
0
\end{bmatrix} .
\end{equation}
Here we have made use of the fact that $\eta$ is the Green's function of a linear differential operator, and that for an $\alpha$-function this takes the form $(1+\alpha^{-1} \d / \d t)^2$.
The network description is completed with the rule for handling jumps
in the synaptic conductance according to $x_i \rightarrow g(x_i)= (\theta_i,s_i, u_i + \alpha)$ whenever a firing event occurs. The times of these events,
$T_i^n$, can be determined from the zeros of an \index{indicator function} indicator function $h :
\RSet^3 \to \RSet$ specified by $h(x_i(T_i^n)) = \theta_i  \! \! \mod 2 \pi - 2 \pi$.  The periodic network state is given by $(\theta_i(t),s_i(t),u_i(t)) = (\theta(t),s(t),u(t)) \equiv x(t)$ for all $i$, with a period $T$ defined by (\ref{Tlinear}).  Here $\theta(t) = \gamma \Gamma \int_0^t P(s) \d s - \Theta t$, $s(t) = P(t)$, and $u(t) = (1+\alpha^{-1} \d /\d t) P(t)$.
By introducing the matrix $M \in \CSet^{3 \times 3}$:
\begin{equation}
M(\beta) = K(T) \e^{[A + \beta \D F] T},
\end{equation}
the MSF for the Lighthouse model can be calculated as
\begin{equation}
\text{MSF} (\beta) = \frac{1}{T} \text{Re} \, \ln (m(\beta)) ,
\end{equation}
where $m(\beta)$ is the largest eigenvalue of $M(\beta)$.  Here, $K(T)$ is a \textit{saltation matrix} that accounts for the jumps in the synaptic perturbations at firing events.  The synchronous state  is stable if the MSF is negative at $\beta=\gamma \widehat{w}_\mu$ where $\mu$ indexes all the eigenvalues of $w$ excluding the one with value $\Gamma$.  Here, $\D F$ is easily calculated whilst the saltation matrix can be computed from the prescription in Appendix A, yielding
\begin{equation}
\D F = \begin{bmatrix}
0 & 1 & 0 \\
0 & 0 & 0 \\
0 & 0 & 0
\end{bmatrix} ,
\qquad 
K(T) = 
\begin{bmatrix}
1 & 0 & 0 \\
\alpha^2/\dot{\theta}(T) & 1 & 0 \\
-\alpha^2/\dot{\theta}(T) & 0 & 1
\end{bmatrix}  ,
\end{equation}
with $\dot{\theta}(T) = \gamma \Gamma P(T) -\Theta$.
\begin{figure}[htbp]
\centering
\includegraphics[width=0.6\textwidth]{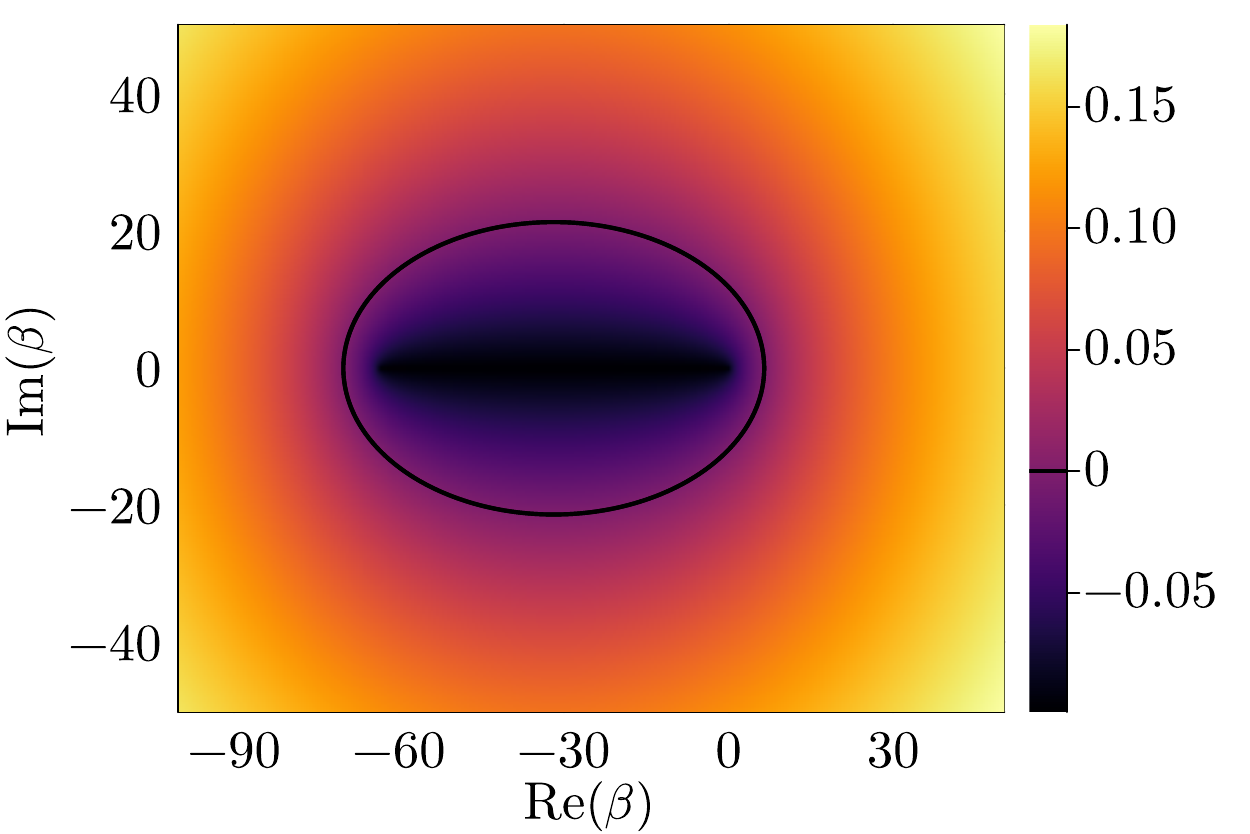}
\caption{The Master Stability Function $\text{MSF} (\beta)$ for a Haken Lighthouse network with $\alpha$-function synapses, shown as a colour plot.
The parameters are $\Theta=-1$, and $\alpha = 0.1$.  The synchronous solution is stable provided all of the eigenvalues of $\gamma w$ lie within the black closed curve.
\label{Fig:MSF}
}
\end{figure}

If the constant value $\Gamma$ is $0$, then we shall say that the system is \textit{balanced}. This scenario is ubiquitous in the modelling of many spiking neural networks \cite{VanVreeswijk1996}.  For a balanced network the existence of the synchronous network state is independent of any of the parameters describing the synaptic interaction.  This means that the form of synaptic coupling cannot induce any nonsmooth bifurcations, such as grazing, since a trajectory $\theta(t) = -\Theta t$ (with period $T= - 2 \pi/\Theta$, $\Theta <0$) cannot tangentially intersect $2 \pi$.  An example MSF is shown in Fig.~\ref{Fig:MSF}, where synchronous solutions are stable in the interior of the region enclosed by the black line.

\begin{figure}
\centering
\includegraphics[width=0.6\textwidth]{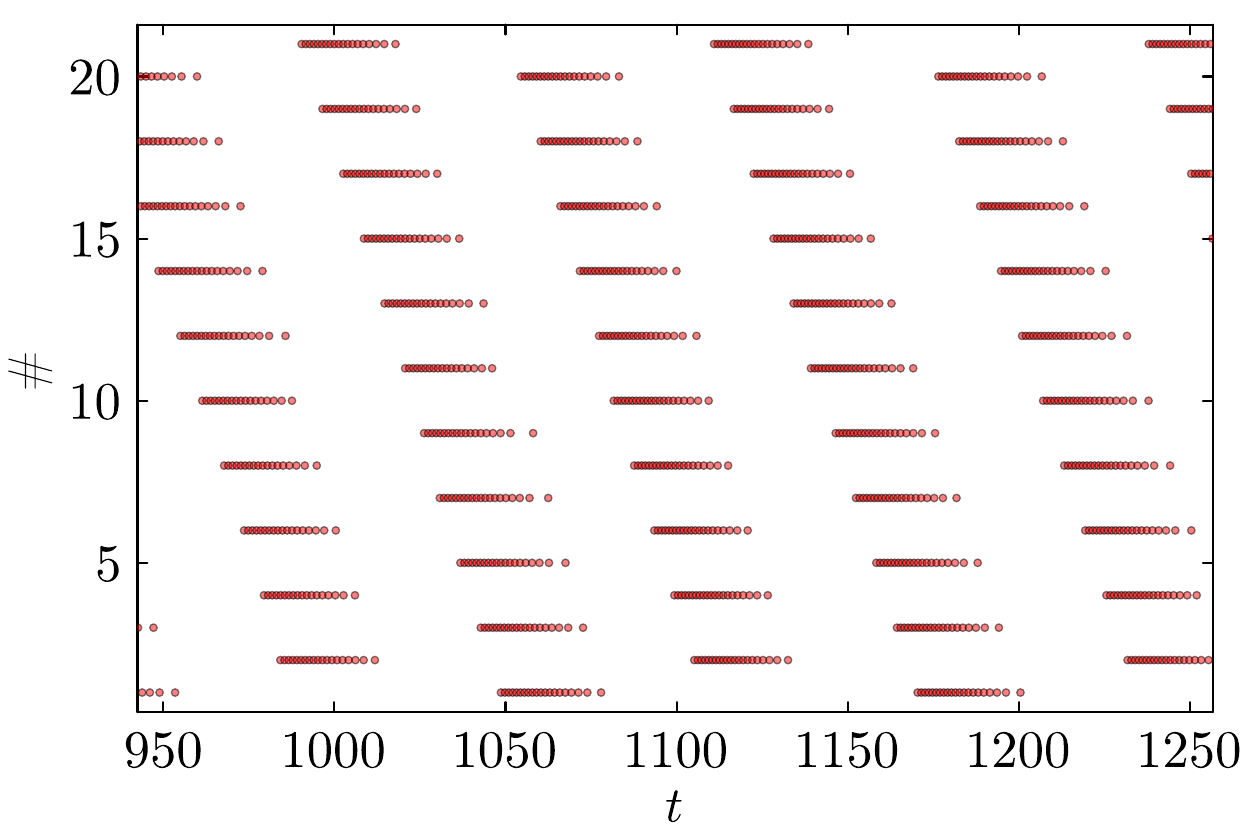}
\caption{A novel bursting behaviour may arise as one leaves the negative region of the MSF along the imaginary axis.  Here we show a simulation of a Haken Lighthouse model with an anti-symmetric circulant connectivity matrix generated by the row vector $(0, -\epsilon, \epsilon, -\epsilon, \epsilon, \ldots, -\epsilon, \epsilon)$ for $N=21$.
The parameters are $\Theta=-1$, $\gamma=1$, $\alpha = 0.1$, and $\epsilon=1.5$ (for which $\beta \simeq 20 \i $).  
\label{Fig:Bursting}
}
\end{figure}

The power of the MSF approach is its ability to determine the stability of the synchronous state for an \textit{arbitrary} network (provided the synchronous solution exists).
When complex values $\gamma \widehat{w}_\mu$ lie within the region where the MSF is negative, then small perturbations to synchronous initial data die away and the system will settle to a synchronous periodic orbit.  When $\gamma \widehat{w}_\mu$ crosses the zero level set of the MSF (from negative to positive) the synchronous solution will destabilise to another firing pattern.  For example choosing an anti-symmetric circulant (Toeplitz) matrix (with $N>2$), which will have purely imaginary (conjugate pairs of) eigenvalues, an instability will occur as one moves from the origin along the imaginary axis (for some sufficiently large value of $\text{Im} \, (\beta)$).  
Circulant matrices have normalised eigenvectors given by $\vec{e}_l = (1,\varpi_l, \varpi_l^2, \ldots, \varpi_l^{N-1})/\sqrt{N}$, where $l=0,\ldots, N-1$, and $\varpi_l=\exp(2 \pi i l /N)$ are the $N$th roots of unity.  The eigenvalues of such matrices are given by $\chi_l = \sum_{j=0}^{N-1} c(|j|) \varpi_l^j$, where $c=(c_0, c_1, \ldots, c_{N-1})$ is the generator (with other rows being obtained by cyclic permutation).
We note that the balance condition enforces $\chi_0=0$.  We also have that $\chi_{N-l}=\chi_l$ for $l=1,\ldots, (N-1)/2$, so that any excited pattern is given by a combination $\vec{e}_m+\vec{e}_{-m} = 2 \text{Re} \, \vec{e}_m$ for some $1\leq m \leq (N-1)/2$.  
In Fig.~\ref{Fig:Bursting} we show a possible bursting wave type of behaviour that can arise in a balanced network with row generator $c = (0, -\epsilon, \epsilon, -\epsilon, \epsilon, \ldots, -\epsilon, \epsilon)$.  For a further discussion of the use of the MSF in studying spiking IF networks, see \cite{Nicks2018}.

\section{A continuum model\label{Sec:Continuum}}

Instead of a discrete labelling of nodes we now turn to a continuous labelling and introduce $\theta = \theta(x,t)$ where $x \in \RSet$ and $t>0$.  The natural extension of the Lighthouse model to the continuum is 
\begin{equation}
\PD{}{t} \theta (x,t) = S(\psi(x,t)),
\label{LH_continuum}
\end{equation}
where
\begin{equation}
\psi(x,t) = \int_{\RSet} w(x,y) s(y,t-\tau(x,y)) \d y, \qquad s(x,t) = \sum_{m \in \ZSet} \eta(t-T^m(x) .
\label{psi_integral}
\end{equation}
Here, $w(x,y)$ specifies the strength of interaction between neuron densities at $x$ and $y$, with $\tau(x,y)$ representing the associated axonal communication delay.
By choosing $\eta$ as the Green's function of a linear differential operator $Q$ such that $Q \eta =\delta$, we may rewrite $\psi$ as an integro-differential equation of the form 
\begin{equation}
Q \psi(x,t) =  \int_{\RSet} w(x,y) \sum_{m \in \ZSet}\delta (t-T^m(y)- \tau(x,y)) \d y .
\label{psi_differential}
\end{equation}
Hence, the continuum Lighthouse model is defined by (\ref{LH_continuum}) using either (\ref{psi_integral}) or (\ref{psi_differential}).  The latter is useful for making a link to more well established firing rate neural fields.  To see this consider the choice of an $\alpha$-function synapse model and restrict attention to small $\alpha$, namely a slow synapse.  In this case a short term temporal average of (\ref{psi_differential}), denoted by $\left \langle \cdot \right \rangle_t = \Delta^{-1} \int_{t}^{t+\Delta} \cdot \, \d s$ for some window $\Delta$ yields
\begin{align}
\left \langle  Q \psi(x,t) \right \rangle_t  =  \int_{\RSet} w(x,y) \left \langle \sum_{m \in \ZSet}\delta (t-T^m(y)- \tau(x,y)) \right \rangle_t \d y .
\label{}
\end{align}
From this we obtain the approximation
\begin{align}
Q \psi(x,t) \simeq  \int_{\RSet} w(x,y) R (y, t - \tau(x,y)) \d y ,
\label{}
\end{align}
where $R(x,t)$ counts the number of spikes per time $\Delta$ at position $x$ at time $t$.  Namely $R$ is a firing rate.  From equation (\ref{slow_period}) we see that a natural choice for this rate function is of the form $\dot{\theta}/(2 \pi)$, which gives the rate based description of the Lighthouse model for slow synapses as
\begin{equation}
Q \psi(x,t) =  \frac{1}{2 \pi} \int_{\RSet} w(x,y) S (\psi(y,t - \tau(x,y)) \d y .
\label{RateField}
\end{equation}
This is in the form of a standard neural field equation and as such may be analysed with a variety of existing techniques, such as those reviewed in \cite{Coombes2014}.  For example, with $w(x,y) = w(|x-y|)$ (distant dependent interactions) and $\tau(x,y) = \tau$ (single fixed delay) a homogeneous steady state (if it exists) will satisfy
$\overline{\psi} = S(\overline{\psi}) \widehat{w}(0)/(2 \pi)$.  Here, anticipating the need for Fourier transforms in a linear stability analysis, the value of $\int_{\RSet} w(y) \d y$ has been expressed as $\widehat{w}(0)$, where
\begin{equation}
\widehat{w}(k) = \int_{-\infty}^\infty \d x w(x) \e^{-\i kx} . 
\end{equation}
A linear stability analysis with respect to perturbations of the form $\e^{\lambda t} \e^{\i k x}$ for $\lambda \in \CSet$ and $k \in \RSet$ yields a spectral equation
\begin{equation}
\left ( 1 + \frac{\lambda}{\alpha} \right )^2 =\frac{S'(\overline{\psi}) \widehat{w}(k)}{2 \pi}\e^{-\lambda \tau}.
\label{lambdaslow1}
\end{equation}
This is recognised as a counterpart to equation (\ref{lambdaslow}).  The latter describes a spectral problem for an homogeneous oscillatory solution of a spiking Lighthouse model, whilst (\ref{lambdaslow1}) describes a spectral problem for an homogeneous steady state of a firing rate Lighthouse model.  Despite one being for a graph and the other a continuum there is no major difference.  The solution $\overline{\psi}$ is linearly stable if $\text{Re} \, \lambda(k) < 0$ for all $k \neq 0$, and the techniques for analysing (\ref{lambdaslow}) may be re-used.  So a real instability is expected when $1 = S'(\overline{\psi})  \max_{k} \widehat{w}(k)/(2 \pi)$ and a dynamic one when
$S'(\overline{\psi}) \min_k \widehat{w}(k)/(2 \pi) -1 \simeq (\pi/(2 \tau)^2$ ($\max_{k} \widehat{w}(k) >0 >\min_{k} \widehat{w}(k)$).

For certain choices of the connectivity function $w$ it is also possible to recast the non-local term defined by (\ref{psi_integral}) (for both the spike and rate model) as the solution to a damped inhomogeneous wave equation, often referred to as the \text{brain wave} equation \cite{Nunez1974}.  This is very useful for performing numerical simulations, though more general approaches for treating the non-local model are also available \cite{Hutt2010}.
 Interestingly, the brain wave equation has also been studied by Haken \cite{Jirsa1997}, though not in the context of the Lighthouse model.  In Appendix B, we present the brain wave equation for $w(x) = \e^{-|x|/\sigma}/(2 \sigma)$ in one spatial dimension and $w(r) =\e^{-r/\sigma}/(2 \pi \sigma^2)$ in two spatial dimensions.  In both cases we note that more general kernels built from linear combinations of exponential functions lead to a set of coupled brain wave equations.
 
\begin{figure}
\centering
\includegraphics[width=0.6\textwidth]{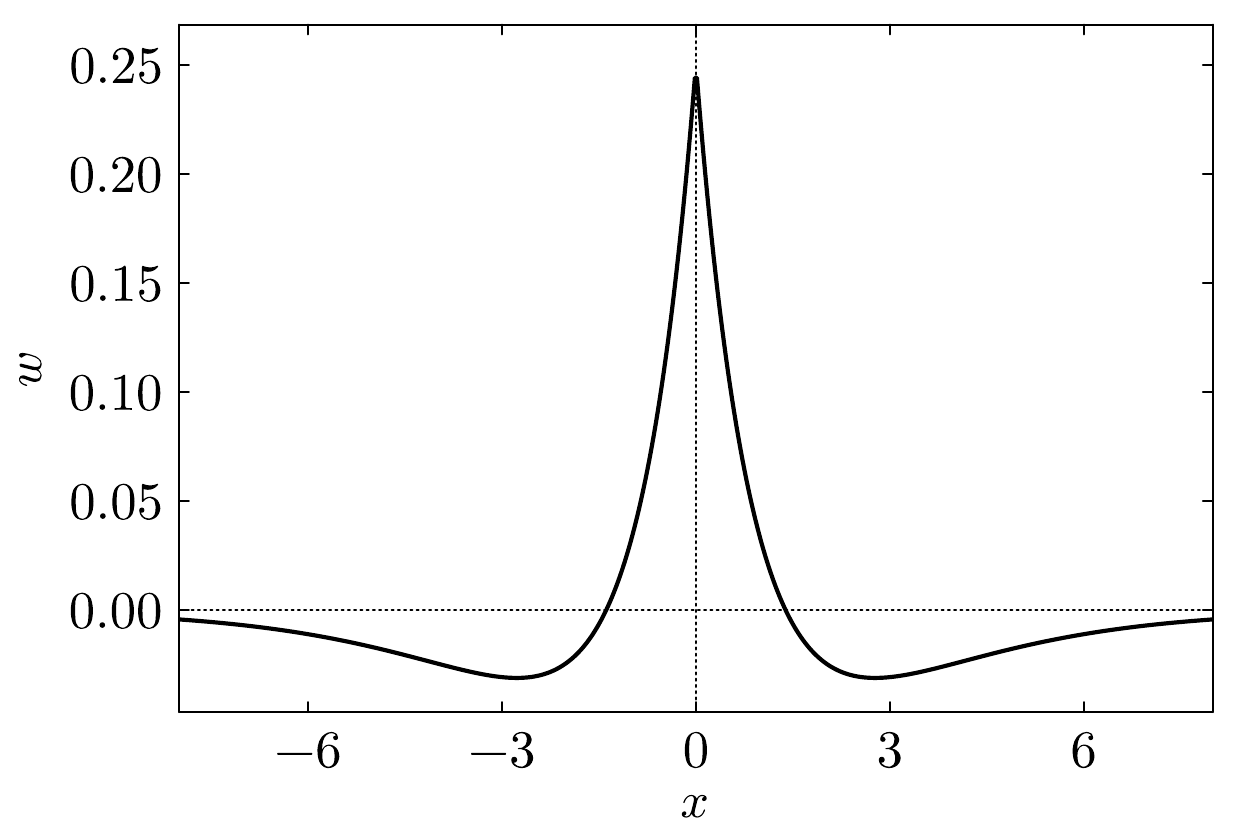}
\caption{An example of an effective Wizard hat connectivity, generated with the choice (\ref{WizardHat}) with $\sigma=2$, $A=1$, and $\Gamma = 0$.
\label{Fig:WizardHat}
}
\end{figure}

\section{Periodic travelling waves\label{Sec:TW}}

Here we consider the spiking continuum model defined by (\ref{LH_continuum}) and (\ref{psi_integral}) and restrict attention to the case of space dependent interactions so that $w(x,y) = w(|x-y|)$ and $\tau(x,y) = |x-y|/v$ for some fixed communication speed $v$.  With the expectation that the model supports periodic travelling waves we look for solutions parameterised by firing times of the form $T^m(x) = mT +  \rho x$.  To determine the dispersion curve $T = T(\rho)$ we substitute this ansatz into (\ref{psi_integral}) to obtain $\psi(x,t) = \psi(\xi)$, where $\xi=t-\rho x$ is a travelling wave frame, and
\begin{equation}
\psi(\xi;T,\rho) = \int_{\RSet} w(|y|) P(\xi - \rho y - |y|/v) \d y , \qquad P(\xi) = \sum_{m \in \ZSet} \eta(\xi -mT) .
\end{equation}
Since $\psi(\xi)$ is $T$-periodic we may also represent it as a Fourier series $\psi(\xi; T,\rho) = \sum_n \psi_n(T,\rho) \e^{2 \pi \i n \xi/T}$, where
\begin{equation}
\psi_n(T,\rho) = \frac{1}{T} \widehat{\eta} (\omega_n) \left [ \widehat{w}_+ (\omega_n(v^{-1}+\rho)) 
+ \widehat{w}_+ (\omega_n(v^{-1}-\rho)) \right ], \qquad \omega_n = 2 \pi n/T ,
\end{equation}
where $\widehat{w}_+(k)$ is a Fourier transform on the half-space:
\begin{equation}
\widehat{w}_+(k) = \int_0^\infty \d y \, w(y) \e^{-\i k y} .
\label{w+}
\end{equation}
A solution for $\theta(x,t)$ can be constructed as $\theta(t - \rho x)$ where
\begin{equation}
\theta (\xi) = \theta(0) + \int_0^\xi S(\psi(s;T,\rho)) \d s.
\label{thetaTW}
\end{equation}
Finally, the dispersion relationship is fixed by the condition $\theta (T) - \theta(0) = 2 \pi$.

\begin{figure}
\centering
\includegraphics[width=0.6\textwidth]{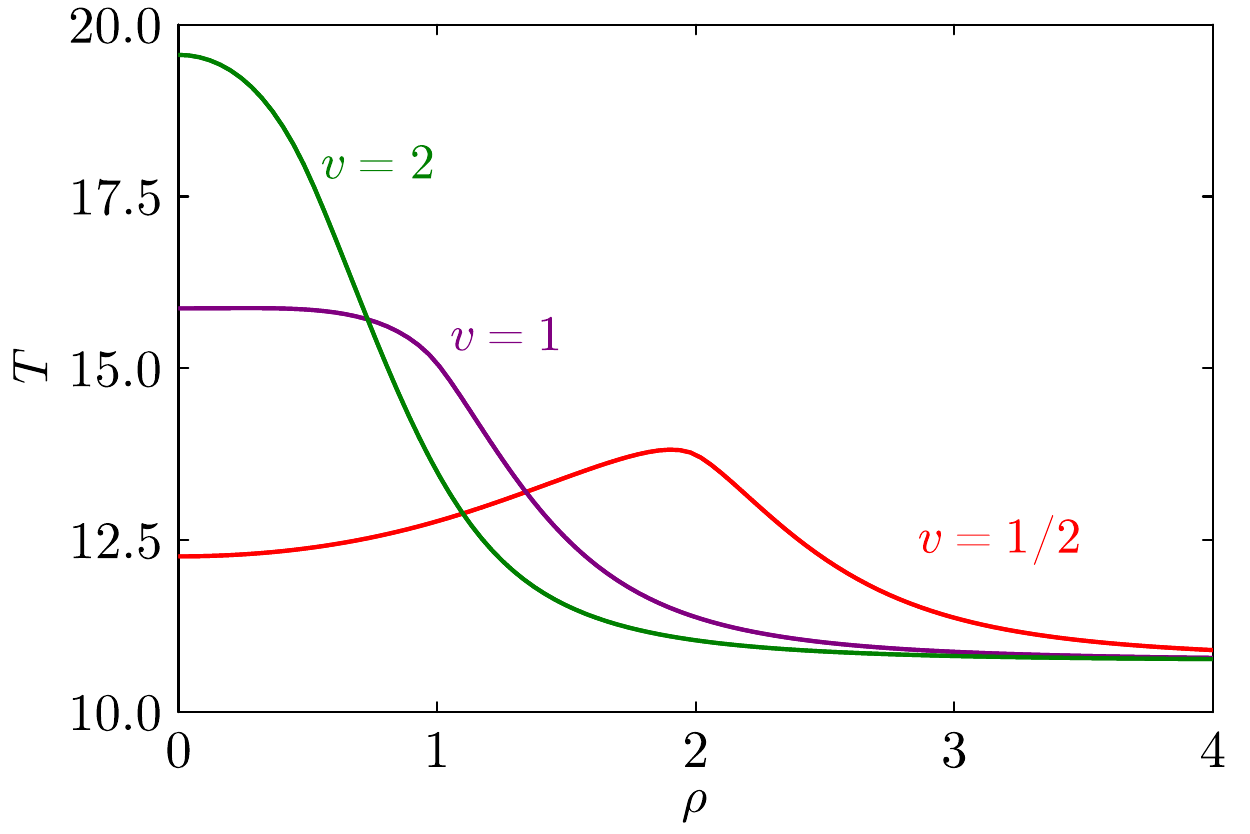}
\caption{A set of dispersion curves, $T=T(\rho)$, for periodic travelling waves in a continuum Haken Lighthouse model with space dependent delays determined by an axonal speed $v$, and with an effective connectivity given by (\ref{WizardHat}).
The parameters are $h=-1$, $r=2$, $\tau=0$, $\Gamma=10$, $\alpha=1$, $\sigma=2$, and $A=1$.
\label{Fig:WaveDispersion}
}
\end{figure}

To describe an effective connectivity function $w(x)$ it is useful to introduce the normalised function $E(x;\sigma) = \e^{-|x|/\sigma}/{2 \sigma}$
and write 
\begin{equation}
\label{WizardHat}
w(x) = A[E(x;1) - (1-\Gamma/A) E(x;\sigma))].
\end{equation} 
We recognise $\Gamma$ as the area under the curve $w(x)$, namely $\Gamma = \widehat{w}(0)$.  For example with $\sigma=2$, $A=1$, and $\Gamma=0$ we can generate a balanced wizard hat shape, such as the one shown in Fig.~\ref{Fig:WizardHat}.  For this choice $\widehat{w}_+(k) = A[\widehat{E}_+(k;1) - (1-\Gamma/A) \widehat{E}_+(k;\sigma))$, with $\widehat{E}_+(k;\sigma) = 1/(2(1+\i \sigma k))$.  Using this in conjunction with the condition $2 \pi - \int_0^T S(\psi(s;T,\rho)) \d s = 0$, we may determine dispersion curves, as illustrated in Fig.~\ref{Fig:WaveDispersion}.

The synchronous state is described by the choice $\rho=0$, and in the next section we treat the linear stability of this solution.

\section{Turing instability analysis of the synchronous state\label{Sec:Turing}}

The study of pattern formation in spiking neural field models is surprisingly little studied, at least without recourse to mean-field reductions \cite{Senk2020}.  Recent work by Dumont and Tarnicerui considers this for renewal neurons, though also uses mean-field equations and not spikes \textit{per se} \cite{Dumont2024}.  However, a Turing instability analysis of the synchronous spiking state in a spatially continuous Lighthouse network can be developed by adapting the techniques for discrete networks described in Sec.~\ref{Sec:SyncandStab}.  Here we note that in contrast to standard Turing instability approaches we consider the linear stability of oscillatory solutions to spatially periodic perturbations, and not that of the homogeneous time-independent steady state.

From (\ref{thetaTW}) a synchronous solution with $\theta(x,t) = \theta(t)$ has an emergent period $T$ that satisfies 
\begin{equation}
2 \pi  = \int_0^T S(\psi(s;T,0)) \d s,
\label{Tsync}
\end{equation}
with the Fourier coefficients of $\psi(t;T,0)$ reducing to $\psi_n(T,0) = \widehat{\eta}(\omega_n) \widehat{w}(\omega_n v^{-1})/T$.
Proceeding in a analogous fashion to Sec.~\ref{sec:linstab}, we introduce a perturbed trajectory $\widetilde{\theta}(x,t)$ that results from perturbations in the firing times.  We denote the perturbed firing times by $\widetilde{T}^m(x)$.  We further introduce the state and time deviations $\delta\theta(x,t) = \widetilde{\theta}(x,t) - \theta(t)$ and $\delta T^m(x) = \widetilde{T}^m(x) - mT$ respectively.  We may relate these state and time deviations by using the firing condition:
\begin{equation}
 \theta (m T)  = \widetilde{\theta}(\widetilde{T}^m(x)) .
 \label{fire-fire1}
\end{equation}
A Taylor expansion of (\ref{fire-fire1}) for small deviations gives the result
\begin{equation}
\delta T^m(x) = - \frac{\delta \theta (x,mT) }{\dot{\theta} (T)} .
\label{dT1}
\end{equation}
Integrating the perturbed trajectory from $\widetilde{T}^m(x)$ to $\widetilde{T}^{m+1}(x)$ gives
\begin{equation}
\widetilde{\theta}(x,\widetilde{T}^{m+1}) - \widetilde{\theta}(x,\widetilde{T}^m) =
\int_{\widetilde{T}^m(x)}^{\widetilde{T}^{m+1}(x)} \d t\,
S \left (
\int_{\RSet} \d y \, w(y) \sum_{p \in \ZSet} \eta(t - \widetilde{T}^p(x-y) - |y|/v)
\right) .
\end{equation}
A Taylor expansion around the synchronous state then gives the linear integro-difference equation
\begin{align}
&\delta \theta(x,(m+1)T) - \delta \theta(x,mT) = \nonumber \\
& \frac{1}{\dot{\theta}(T)} \int_0^T \d t \, S'(\psi(t))
\sum_{p \in \ZSet} \int_{\RSet} \d y \, w(y)  \eta'(t + pT - |y|/v) [\delta \theta(x-y,(m-p)T) - \delta \theta(x,mT)] .
\label{linear_sync}
\end{align}
Noting that periodic functions are eigenfunctions of the convolution operator and that a power function is the natural solution of a linear difference equation,
the continuous spectrum can be found by taking solutions of the form 
$\delta \theta(x,mT) = \e^{\i kx} \e^{m \lambda}$ where $\lambda \in \CSet$ and $k \in\RSet$, and can be represented as the zeros of the complex function $\mathcal{E}(\lambda, k)$, where 
\begin{equation}
\mathcal{E}(\lambda, k) = (\e^{\lambda} - 1) \dot{\theta}(T) - (a(\lambda,k) - a(0,0)).
\label{Econtinuum}
\end{equation}
Here,
\begin{equation}
a(\lambda,k) = \frac{1}{T} \sum_{n \in \ZSet}  f_n(\lambda) \widehat{\eta} (\omega_n - \i \lambda/T) 
\left [ \widehat{w}_+ ((\omega_n -\i \lambda/T) v^{-1} -k) 
+ \widehat{w}_+ ((\omega_n -\i \lambda/T) v^{-1} +k) 
\right ] ,
\end{equation}
with, 
\begin{equation}
f_n(\lambda) = \int_0^T \d t \, S'(\psi(t)) \FD{}{t}\e^{\i [\omega_n -\i \lambda /T]t}  = S'(\psi(T)) (\e^\lambda - 1) -  \int_0^T \d t \, S''(\psi(t)) \dot{\psi}(t) \e^{\i [\omega_n -\i \lambda /T]t} ,
\label{fn}
\end{equation}
and we note that if $|\dot{\psi}(t)| \ll 1$ for all $t$ then $f_n(\lambda) = S'(\psi(T)) (\e^{\lambda}-1)$ for all $n$ (and this is always true for the linearised Lighthouse model since $S''(\psi)=0$ for all $\psi$).

A Turing instability can occur if $\lambda \in \RSet$.  In this case the bifurcation condition is $\mathcal{E}(0, k_c) =0$ for some critical wavenumber $k_c$. This can be found by demanding that $\lambda(k)$ tangentially touches zero from below, or equivalently that $\left. \partial \mathcal{E} (\lambda, k) / \partial k \right |_{\lambda=0} = 0$ with $\text {Re} \, \lambda (k) \leq 0$.
From (\ref{Econtinuum}) we note that for $v^{-1}=0$ (no communication delay) a stationary point occurs when $\widehat{w} ' (k) = 0$. Thus, in this case $k_c$ is determined solely by the properties of the spatial connectivity pattern.
Beyond a Turing instability we expect to excite patterns of the form $\e^{\i k_c x}$.  For $\lambda \in \CSet$ Hopf, and Turing-Hopf bifurcations are possible.  The former requires  $\mathcal{E}(\i \omega_c, 0) =0$ and the latter  $\mathcal{E}(\i \omega_c, k_c) =0$, for some critical frequency $\omega_c$.  In all cases we must further ensure that the spectrum has $\text {Re} \, \lambda (k) \leq 0$ and only touches the imaginary axis from the left.  The condition for the spectrum to tangentially touch the imaginary axis is
\begin{equation}
\left | 
\frac{\partial( \text{Re} \, \mathcal{E}(\i \omega, k), \text{Im} \, \mathcal{E}(\i \omega, k))}{\partial ( w, k)}
\right | = 0 .
\label{tangency}
\end{equation}

In the limit of slow synapses we may make the same arguments as presented in Sec.~\ref{sec:slow} to show that the non-trivial continuous spectrum is defined by
\begin{equation}
\frac{1}{\widehat{\eta}(- \i \lambda /T)} = \frac{S'(\widehat{w}(0)/T)}{2 \pi} \widehat{w}(k -\i \lambda/(vT)),
\end{equation}
where the period of the synchronous orbit is given by $2 \pi/T = S(\widehat{w}(0)/T)$.
The condition for a Turing bifurcation is independent of $v$ and given by $1 = S'(\widehat{w}(0)/T) \max_k \widehat{w}(k) /(2 \pi)$.
The condition for a Hopf bifurcation is a solution to the complex equation $1/{\widehat{\eta}(\omega_c/T)} = S'(\widehat{w}(0)/T)\widehat{w}(\omega_c/(vT))/(2 \pi)$, which may be solved for a pair $(\omega_c,a)$ for some system parameter $a$.  Similarly, for the Turing-Hopf bifurcation we must solve the system $1/{\widehat{\eta}(\omega_c/T)} = S'(\widehat{w}(0)/T)\widehat{w}(k_c+\omega_c/(vT))/(2 \pi)$ for $(\omega_c,k_c,a)$, where we must also make use of (\ref{tangency}).

\begin{figure}[htbp]
\centering
\includegraphics[width=0.6\textwidth]{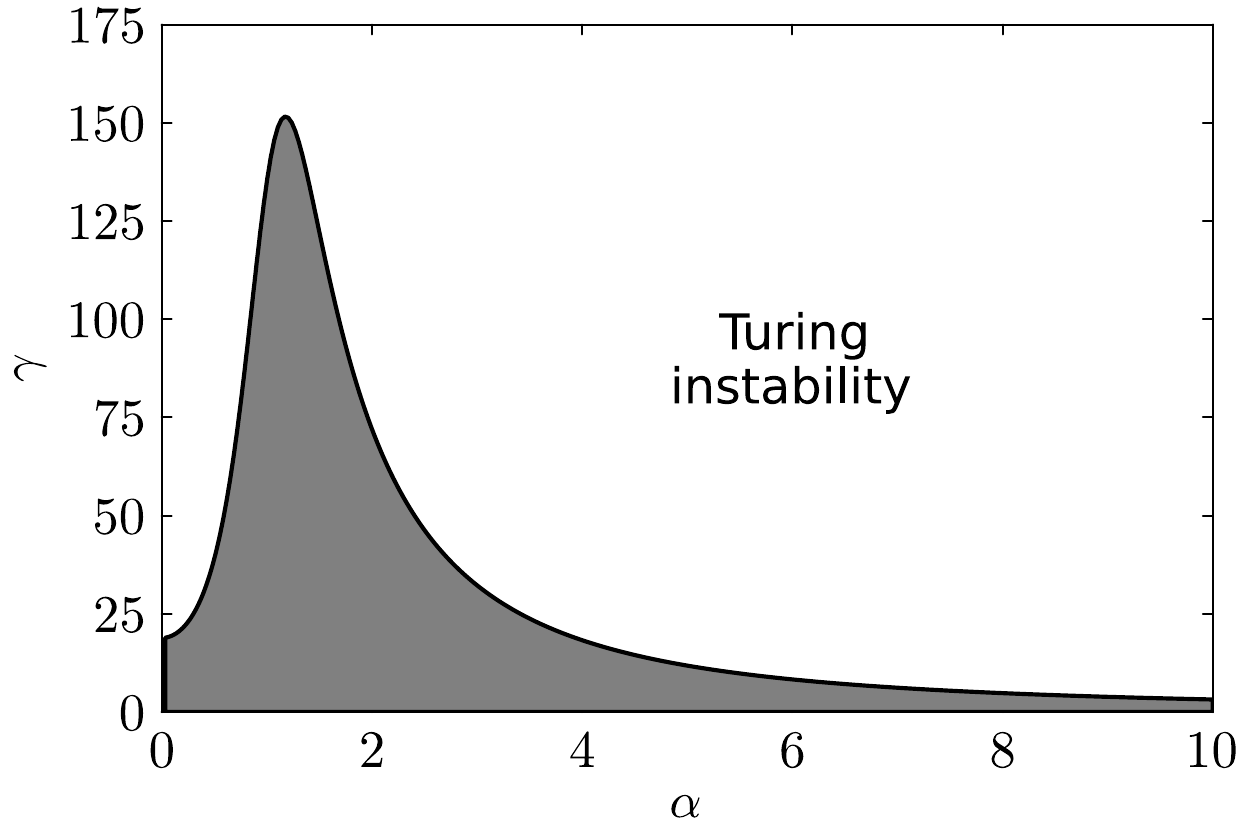}
\caption{The critical curve $\gamma = \gamma_c(\alpha)$ for a Turing bifurcation in a balanced network prescribed by (\ref{WizardHat} with $v^{-1}=0$ and $S=S_L$.
The synchronous spiking network state is stable in the grey zone.
The parameters are $\Theta = -1$, $\sigma=2$, and $A=1$.
\label{Fig:StaticTuring}
}
\end{figure}
\begin{figure}
\centering
\includegraphics[width=0.6\textwidth]{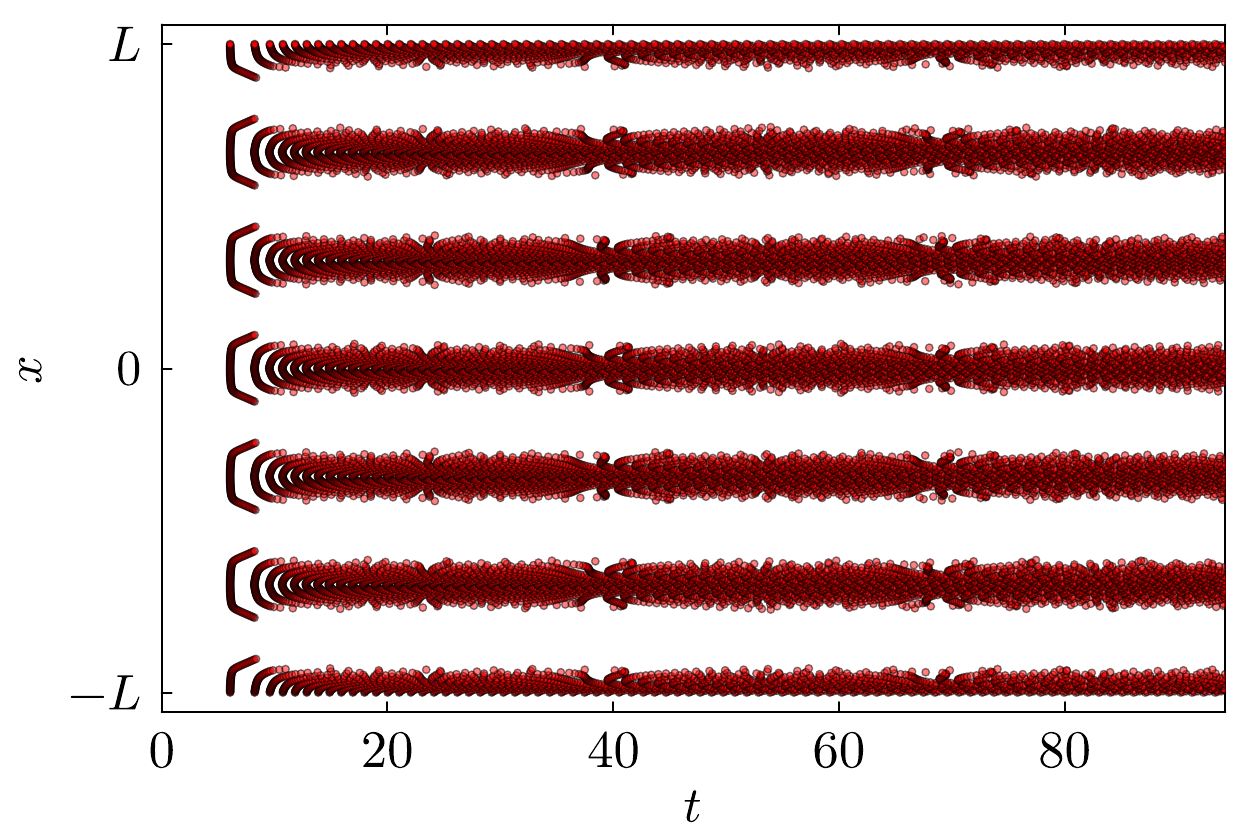}
\caption{A raster plot of a solution that can emerge beyond a Turing bifurcation.  Simulations were performed using a regular spatial discretisation on an interval $2 L = 6 (2\pi/k_c)$ using $2^{10}$ mesh points and periodic boundary conditions.
The parameters are as in Fig.~\ref{Fig:StaticTuring} with $\alpha=4$, and $\gamma=25$. 
\label{Fig:TuringRaster}
}
\end{figure}

By way of illustration, consider a balanced network prescribed by (\ref{WizardHat}) with $v^{-1}=0$ and $S=S_L$.  In this case a Turing bifurcation is defined by the condition
\begin{equation}
-\Theta = \gamma \widehat{w}(k_c) P(T), \qquad k_c = \frac{1}{\sqrt{\sigma}}, \qquad T = - \frac{2 \pi}{\Theta} .
\label{StaticTuring}
\end{equation}
Here, $\widehat{w}(k) = A [\widehat{E}(k,1) - \widehat{E}(k,\sigma)]$, and $\widehat{E}(k,\sigma) = (1+\sigma^2 k^2)^{-1}$.
We may use (\ref{StaticTuring}) to define a critical curve $\gamma = \gamma_c(\alpha)$, above which the synchronous solution will give way to a spatially modulated spiking pattern with a wavelength $\sim 2 \pi \sqrt{\sigma}$ close to bifurcation.  A plot of such a critical curve is shown in Fig.~\ref{Fig:StaticTuring}, and in Fig.~\ref{Fig:TuringRaster} we show a raster plot of a pattern that can emerge beyond the instability.  


\subsection{Extension to two spatial dimensions}

The analysis above is readily generalised to two spatial dimensions where $x \rightarrow \vec{r} \in \RSet^2$.  Choosing a radially symmetric kernel such that 
$w(\vec{r}) = w(r)$ and  $r=|\vec{r}|$, the above analysis generalises under the replacement $\psi_n(T,0) \rightarrow \widehat{\eta}(\omega_n) \widehat{w}(0,\omega_n v^{-1})/T$, where
\begin{equation}
\widehat{w} (k,a) = \int_{\RSet^2} \d \vec {r} \, w(|\vec{r}|) \e^{-\i a |\vec{r}|} \e^{-\i \vec{k} \cdot \vec{r}} ,
\end{equation}
with $k  = |\vec{k}|$ for $\vec{k} \in \RSet^2$.
The period of the synchronous solution in two spatial dimensions is determined by (\ref{Tsync}), and similarly the expressions for determining stability
can be expressed in terms of $\widehat{w}(k,(\omega_n - \i \lambda/T)v^{-1})$.
For example we might introduce $E(r;\sigma) = \e^{-r/\sigma}/(2 \pi \sigma^2)$ and choose $w(r) = A[E(r;1) - (1-\Gamma/A) E(r;\sigma))]$, for which
$\widehat{w} (k,a) = A[ \widehat{E}(k,a;1) - (1-\Gamma/A) \widehat{E}(k,a;\sigma))]$, where
\begin{equation}
\widehat{E}(k,a;\sigma) =\frac{1}{\sigma^2} \frac{\sigma^{-1} +\i a}{((\sigma^{-1} +\i a)^2+k^2)^{3/2}}.
\end{equation}

\subsection{Travelling  wave stability}

The Turing instability analysis for a synchronous state is easily generalised to cover the periodic travelling waves described in Sec.~\ref{Sec:TW}.  This leads to a slight generalisation of (\ref{linear_sync}) to the form
\begin{align}
&\delta \theta(x,(m+1)T) - \delta \theta(x,mT) = \nonumber \\
& \frac{1}{\dot{\theta}(T)} \int_0^T \d t \, S'(\psi(t;T,\rho))
\sum_{p \in \ZSet} \int_{\RSet} \d y \, w(y)  \eta'(t + pT + \rho y- |y|/v) [\delta \theta(x-y,(m-p)T) - \delta \theta(x,mT)] .
\label{linear_TW}
\end{align}
To analyse this it is useful to consider eigenfunctions of the linear integro-difference operator defined by the the first expression on the right hand side of (\ref{linear_TW}) (and see \cite{Bressloff2000a} for a related discussion), namely $\mathcal{A} \ast u = a u$, where $u = (\ldots, u_m(x), \ldots)$ and
\begin{equation}
[\mathcal{A}\ast u] (x,m) = \sum_{p \in \ZSet}
\int_0^T \d t \, S'(\psi(t;T,\rho)) \int_{\RSet} \d y \, w(y)  \eta'(t + pT + \rho y- |y|/v) u_{m-p} (x-y) .
\label{A}
\end{equation}
This has eigenfunctions of the form $u_m(x) = \e^{\i k x} \e^{\lambda [m + \rho x/T]}$ with corresponding eigenvalues $a(\lambda,k;\rho)$:
\begin{equation}
\frac{1}{T} \sum_{n \in \ZSet}  f_n(\lambda) \widehat{\eta} (\omega_n - \i \lambda/T) 
\left [ \widehat{w}_+ ((\omega_n -\i \lambda/T) v^{-1} -k + \rho \omega_n) 
+ \widehat{w}_+ ((\omega_n -\i \lambda/T) v^{-1} +k - \rho \omega_n) 
\right ],
\label{arho}
\end{equation}  
where $f_n(\lambda)$ takes the form (\ref{fn}) and $\widehat{w}_+$ is defined by (\ref{w+}).
For a proof of the above see Appendix C.
Thus, the spectrum can be represented as the zeros of the complex function $\mathcal{E}(\lambda, k; \rho)$, where 
\begin{equation}
\mathcal{E}(\lambda, k; \rho) = (\e^{\lambda} - 1) \dot{\theta}(T) - [a(\lambda,k;\rho) - a(0,0;\rho)] .
\end{equation}

\section{Localised bump solutions\label{Sec:Bump}}

The investigation of localised activity states, often referred to as \textit{bumps}, in fully spiking neural network models with fast synaptic interactions has largely relied on numerical simulations. For instance, the study of IF networks by Laing and Chow revealed several dynamic behaviours that are not captured by traditional firing rate models \cite{Laing2001}. Their results showed that (i) stable bump solutions can emerge in spiking networks (coexisting with the quiescent state) provided the neurons within the bump fire asynchronously, and (ii) as the speed of synaptic processing increases, these bumps can lose stability, giving rise to either wandering bumps or travelling waves.  Interestingly, this dynamic behaviour is not easily seen in rate-based continuum neural field models without the inclusion of stochasticity \cite{Kilpatrick2013,Avitabile2017}.  Similarly, the deterministic Lighthouse model can support localised bump solutions that may undergo a transition to wandering \cite{Chow2006}.

\subsection{Rate-based Amari reduction}

Before considering the spike-based model let us first set the scene with an analysis of localised states in the rate-based model (\ref{RateField}).  Here, we shall focus on the \textit{Heaviside} case, which since the seminal work of Amari has been shown to be particularly amenable for the analysis of spatially localised states \cite{Amari1977,Amari2014}.
We shall consider the continuum model (\ref{RateField}) in the absence of a delay:
\begin{equation}
Q \psi(x,t) =  \frac{1}{2 \pi} \int_{\RSet} w(|y|) S (\psi(x-y,t)) \d y .
\label{Amaristyle}
\end{equation}
Time-independent solutions $\psi(x,t) = q(x)$ satisfy the nonlinear integral equation
\begin{equation} 
q (x) =  \frac{1}{2 \pi} \int_{\RSet} w(|y|) S(q(x-y)) \d y .
\label{Rmodel}
\end{equation}
The nonlinearity in (\ref{Rmodel}) makes it very hard to obtain analytical results for the solution $q$.  However, in the case that $S$ is a Heaviside ($r \rightarrow 0$ in (\ref{S})) results for spatially localised states are readily obtained.   For example, for certain conditions on the connectivity, it is natural to consider a spatially localised solution 
of the form $q(x) \geq h$ for $x_1 <x < x_2$ with $q'(x_1) >0$ and $q'(x_2)<0$ with some decay toward a sub-threshold constant profile as $x \rightarrow \pm \infty$. In this case the exact solution is given simply by
\begin{equation} 
q(x) = \frac{1}{2 \pi}\int_{x_1-x}^{x_2-x} w(|y|) \d y ,
\end{equation} 
with derivative $2 \pi q'(x) = w(|x_1-x|) - w(|x_2-x|)$.
The conditions $q(x_1)=h$ and $q(x_2)=h$ both lead to the equation
\begin{equation}
 h  = \frac{1}{2 \pi}\int_0^\Delta w(y) \d y ,
\label{Delta}
\end{equation} 
describing a family of solutions with $\Delta=(x_2-x_1)$, assuming $w(0)> w(\Delta)$ (to ensure the correct sign of the gradients at $x=x_{1,2}$). 
The stability of this solution is determined by the condition $\text{Re} \, \lambda <0$ where $\lambda \in \CSet$ is a zero of the complex function $\mathcal{E}(\lambda)$ \cite{Coombes2003}:
\begin{equation}
\mathcal{E}(\lambda) = \frac{1}{\widehat{\eta}(- \i \lambda)} - \frac{w(0) \pm w(\Delta)}{w(0) - w(\Delta)} .
\label{ERateBump}
\end{equation}
We note that $\lambda = 0$ is always a solution, reflecting the fact there is a family of bump solutions related by spatial translation.
For $\lambda \in \RSet$ we see that another zero eigenvalue is possible when $w(\Delta) = 0$.  Moreover, from (\ref{Delta}) we see that $2 \pi \, \d h /\d \Delta = w(\Delta)$, so that a bifurcation is expected when $h'(\Delta) = 0$, with stable (unstable) solutions for $w(\Delta)<0$ ($w(\Delta)>0$).  No dynamic instabilities are possible in the rate model.

\subsection{Bumps as waves}

In contrast to IF networks where bumps have been observed numerically \cite{Laing2001,Bressloff2000b} the Lighthouse model is sufficiently simple that bumps can be constructed analytically when the nonlinearity of the model is chosen to be a Heaviside.  For this choice, the main reason for the tractability over the IF model is that all cells within the bump fire with a uniform frequency, whereas in IF networks the bump edges have a highly variable (and non-periodic) firing pattern.

Here we restrict attention to the case of a Heaviside form of nonlinearity $S$ though without restriction on the speed of the synapse.
It will be convenient to work in the integral representation with
\begin{equation}
\psi(x,t) = \int_{\RSet} \d y w(|x-y|) \sum_{m \in \ZSet} \eta(t-T^m(y)) .
\label{psi}
\end{equation}
We focus on the class of (symmetric) maximally firing one-bump solutions parameterised by firing times 
\begin{equation}
T^m(x) = 2 \pi m + \rho |x| ,
\label{ansatz}
\end{equation}
such that $\psi(x,t) \geq h$ for $x \in [x_1, x_2]$ for all $t$ ($H=1$), and
$\psi(x,t) < h$ otherwise ($H=0$).  As in Sec.~\ref{Sec:TW}, the firing ansatz (\ref{ansatz}) is in the form of a travelling wave, albeit this time in a symmetric form that spreads left and right from some origin.
Here we take $x_1 = -\Delta/2 = -x_2$.  We are free to choose such a coordinate system by translational invariance of the continuum model. In this case $\psi(x,t)$ is given by
\begin{equation}
\psi(x,t) = \int_{x_1}^{x_2} \d y \, w(|x-y|) P(t-\rho |y|) , \qquad P(t) = \sum_{m \in \ZSet} \eta(t-2 \pi m),
\end{equation}
subject to the condition
\begin{equation}
\min_t \psi(\Delta/2 , t) = h.
\label{bumpwidth}
\end{equation}
Since, $P(t)$ is $2 \pi$-periodic then $\psi(x,t)$ is also $2 \pi$-periodic and we may write this in the form $\psi(x,t) = \sum_{n \in \ZSet} \psi_n(x) \e^{\i n t}$, with
$\psi_n(x) = P_n F(x,\rho n)$, where $P_n = \widehat{\eta}(n)/(2 \pi)$ and
\begin{equation}
F(x,k) = \int_{x_1}^{x_2} \d y \, w(|x-y|) \e^{-\i k |y|} .
\label{F}
\end{equation}
The bump width $\Delta = \Delta(\rho)$ is then determined as 
\begin{equation}
h =  \frac{1}{2 \pi}\int_{0}^{\Delta} w(y) \d y  + \frac{1}{\pi} \text{Re} \left ( \sum_{n > 0}
\widehat{\eta}(n) G(\rho n) \e^{\i n t^*} \right )  ,
\label{bumpwidth1}
\end{equation}
where $G(k)=F(\Delta/2,k)$, and $t^*$ denotes the value of $t$ for which $\psi(\Delta/2,t)$ has a minimum.
This can be determined from the condition $\sum_{n \in \ZSet} \psi_n(\Delta/2) (\i n) \e^{\i n t^*} = 0$.
We note that for slow synapses equation (\ref{bumpwidth1}) reduces to (\ref{Delta}), as expected.
Thus, we may regard the second term on the right hand side of (\ref{bumpwidth1}) as a correction to the standard Amari firing rate model description that takes into account the
effects of non-slow synaptic processing. 
Note that for $E(x;\sigma)=\e^{-|x|/\sigma}/(2 \sigma)$ we have that
\begin{equation}
\int_{-\Delta/2}^{\Delta/2} \d y \, E(|\Delta/2-y|;\sigma) \e^{-\i k |y|}  = \frac{1}{2 \sigma} \left \{
\frac{\e^{-\Delta/(2 \sigma)} - \e^{-\Delta/\sigma} \e^{-\i k \Delta/2}}{\i k +\sigma^{-1}}
- \frac{\e^{-\i k \Delta/2} - \e^{-\Delta/(2 \sigma)}}{\i k -\sigma^{-1}}
\right \} ,
\end{equation}
so that $G(k)$ is easily calculated for the case $w(x) = A[E(x;1) - (1-\Gamma/A) E(x;\sigma)]$.

In terms of the phase variables we may write the localised bump solution as
\begin{equation}
\PD{}{t} \theta (x,t) = \begin{cases}
1 & x \in [-\Delta/2 , \Delta/2 ] \\
0 & \text{otherwise}
\end{cases}, 
\qquad
\theta (x,t) = \begin{cases}
t - \rho |x| & x \in [-\Delta/2 , \Delta/2 ] \\
c(x) & \text{otherwise}
\end{cases} ,
\label{bump}
\end{equation}
for some arbitrary function $c(x)$.

\begin{figure}
\centering
\includegraphics[width=0.4\textwidth]{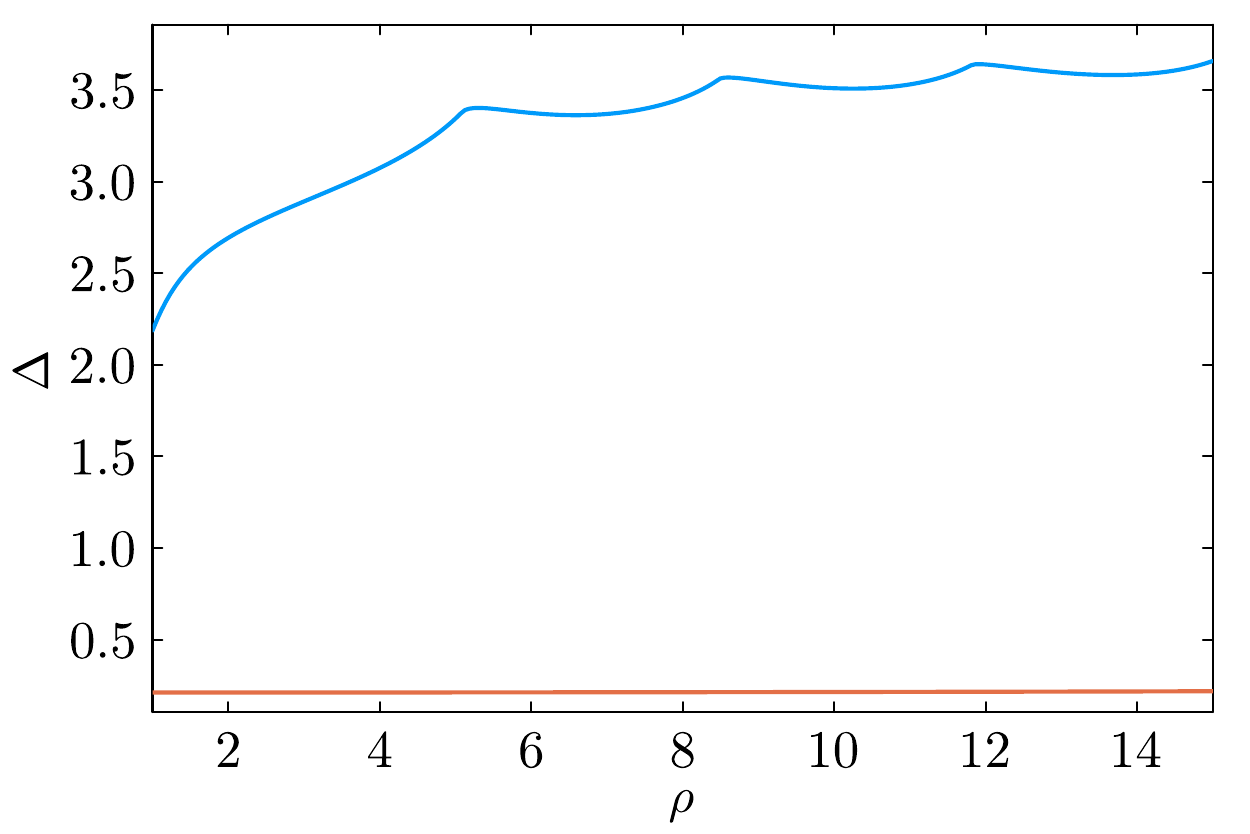} \hspace*{1cm}
\includegraphics[width=0.4\textwidth]{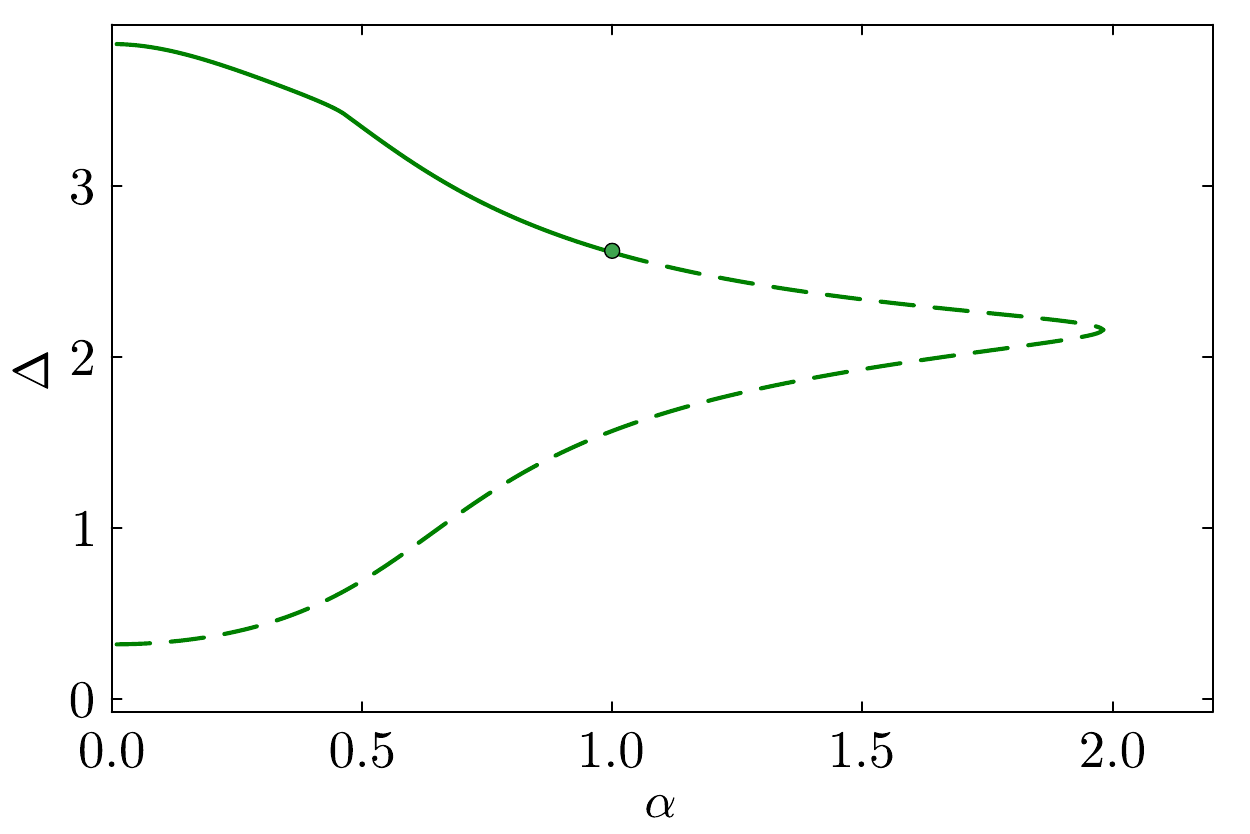}
\caption{
\label{Fig:Deltavsrho}
A plot of the theoretical predictions for bump width, with parameters:  $h=0.01$, $A=1$, $\Gamma=0$, $\sigma=2$, and $M=50$.
Left: Solution branches for $\Delta=\Delta(\rho)$ with $\alpha = 0.5$. An examination of the spectra from Fig.~\ref{Fig:BumpSpectrum} shows that the upper branch (blue curve) is stable and the lower branch (red curve) is unstable at $\rho=5$.
Right: Solution branches for the bump solution with $\rho=5$ showing $\Delta=\Delta(\alpha)$.
}
\end{figure}
In Fig.~\ref{Fig:Deltavsrho} we show a plot of the theoretical bump width $\Delta = \Delta(\rho)$ for fixed $\alpha$ and also a plot of $\Delta = \Delta(\alpha)$ for fixed $\rho$, as determined from (\ref{bumpwidth1}).  As in the Amari rate model we see two branches of solutions.  In the next section we shall treat the stability of these branches.


\subsection{Bump stability}

To determine the linear stability of the spiking bump solution it is convenient to recast (\ref{psi}) in the equivalent form
\begin{equation}
\psi(x,t) = \int_{\RSet} \d y w(|x-y|) \int_0^\infty \d s \, \eta(s) \sum_{m \in \ZSet} \delta(t-s-T^m(y)) .
\label{psi1}
\end{equation}
To relate this back to the dynamics of the phase variable (that generate the spike times) we then use the properties of Dirac-delta functions to obtain the representation
\begin{equation}
\sum_{m \in \ZSet} \delta(t-T^m(x)) = \sum_{m \in \ZSet} |\theta_t(x,t)| \delta(\theta(x,t) - 2 \pi m) . 
\label{Dirac}
\end{equation}
Here, the subscript notation $\theta_t$ is used to denote the partial derivative of $\theta$ with respect to $t$.
The form of (\ref{Dirac}) was much favoured by Haken in his analysis of the Lighthouse model.
A perturbed trajectory $\widetilde{\psi}(x,t)$ may be obtained from (\ref{psi1}) under the replacement $T^m(x) \rightarrow \widetilde{T}^m(x)$ for some perturbed spike times $\widetilde{T}^m(x)$.  Introducing $\delta \psi(x,t) = \widetilde{\psi}(x,t) - \psi(x,t)$, and assuming $\delta \psi(x,t)>0$ we find to first order that
\begin{align}
\delta \psi(x,t) &= \int_\RSet \d y \, w(|x-y|) \int_0^\infty \d s \, \eta(s) \times \nonumber \\
& \sum_{m \in \ZSet} \left \{
\theta_{t} (y, t-s) \delta'(\theta(y,t-s)- 2 \pi m) 
+ \partial_{t} \delta \theta (y, t-s) \delta (\theta(y, t-s) - 2 \pi m)
\right \} .
\label{deltapsi}
\end{align}
Here, $\partial_t \delta \theta (x, t) = S'(\psi(x,t)) \delta \psi(x,t)$.  The first term on the right hand side of 
(\ref{deltapsi}) is directly proportional to $\psi_t(x,t)$, and in the analysis that follows we will see that we only have to consider perturbations at the bump edges.  With this in mind we may further fix $t=t^*$ so that $\psi_t(\pm \Delta/2,t) = 0$ and conveniently restrict attention to only the second contribution from the right hand side of (\ref{deltapsi}).  Assuming solutions of the form $\delta \psi(x,t) = \delta \psi(x) \e^{\lambda t}$ for some $\lambda \in \CSet$, (\ref{deltapsi}) reduces to (before setting $t=t^*$)
\begin{equation}
\delta \psi(x) = \int_\RSet \d y \, w(|x-y|) \int_0^\infty \d s \, \eta(s) 
\sum_{m \in \ZSet} \delta (\psi(y,t-s) - h) \delta \psi(y) \e^{-\lambda s}
\delta (t-s - \rho|y| - 2 \pi m) 
\label{deltapsi1},
\end{equation}
where we have substituted for $\theta$ using (\ref{bump}), and used the result that $S'(x)=H'(x-h) = \delta(x-h)$.
We may now collapse the integral over $y$ to give
\begin{equation}
\delta \psi(x) = \sum_{z \in \{ x_1, x_2 \}} \delta \psi(z) w(x-z) \sum_{m \in \ZSet} \int_0^\infty \d s \, \eta(s) \e^{-\lambda s} \frac{\delta (t^*-s-\rho |z| - 2 \pi m)}{|\psi_x(z,t^*)|}.
\label{deltapsi2}
\end{equation}
Next we collapse the integral over $s$ to give
\begin{equation}
\delta \psi(x) = \sum_{z \in \{ x_1, x_2 \}} \frac{\delta \psi(z) w(|x-z|)}{|\psi_x(z,t^*)|} \sum_{m \in \ZSet} \eta(t^*-\rho |z| - 2 \pi m) \e^{-\lambda (t^*-\rho |z| - 2 \pi m)} .
\label{deltapsi3}
\end{equation}
To perform the summation over $m$ in (\ref{deltapsi3}) we use the integral representation 
$\eta(t) = (2 \pi)^{-1} \int_{\RSet}\d k \, \widehat{\eta}(k) \e^{\i k t}$ and make use of a Dirac comb in the form 
$\sum_m \e^{2 \pi \i m [k- \i \lambda]} = \sum_p\delta (k -\i \lambda + p)$.
Thus, (\ref{deltapsi3}) simplifies to 
\begin{equation}
\delta \psi(x) = \sum_{z \in \{ x_1, x_2 \}} \frac{\delta \psi(z) w(|x-z|)}{2 \pi |\psi_x(z,t^*)|} \sum_{p \in \ZSet} \widehat{\eta} (p - \i \lambda) \e^{\i p [t^*-\rho|z|]},
\label{deltapsi4}
\end{equation}
Setting $x=x_{1,2}$ in (\ref{deltapsi4}) gives a pair of linear equations for $(\delta \psi(x_1),\delta \psi(x_2)$ that will only have a non-trivial solution if $\det (\mathcal{A}(\lambda) - I_2)=0$, where
\begin{equation}
\mathcal{A} (\lambda) = \frac{\sum_{p \in \ZSet } \widehat{\eta} (p - \i \lambda) \e^{\i p [t^*-\rho \Delta/2]}}{2 \pi| \psi_x(\Delta/2, t^*)|}
\begin{bmatrix}
w(0) & w(\Delta) \\
w(\Delta) & w(0)
\end{bmatrix}.
\end{equation}
Thus, the spectrum is defined by the two equations $\mathcal{E}_\pm (\lambda) = 0$, where
\begin{equation}
\mathcal{E}_\pm (\lambda) = 1 - \frac{1}{2 \pi P(t^*-\rho \Delta/2)} \frac{w(0) \pm w(\Delta) }{|w(0) - w(\Delta)|} \sum_{p \in \ZSet } \widehat{\eta} (p - \i \lambda) \e^{\i p [t^*-\rho \Delta/2]},
\label{ESpikeBump}
\end{equation}
and we have used the result that $\psi_x(x,t) = w(x_1-x) P(t - \rho |x_1|) -w(x_2-x)P(t-\rho|x_2|) + \psi_t(x,t)$.  Assuming $w(0)>w(\Delta)$ we see that $\mathcal{E}_-(0)=0$.  Thus, there is always an eigenvalue with $\lambda=0$ as expected from translation invariance.
We note that for slow synapses (where we only need retain the $p=0$ term and $P(t) = (2 \pi)^{-1}$ for all $t$), the spectrum defined by (\ref{ESpikeBump}) reduces to that described by (\ref{ERateBump}) as expected.  Unlike in the rate model it is possible for dynamic instabilities to appear.

\begin{figure}
\centering
\includegraphics[width=0.4\textwidth]{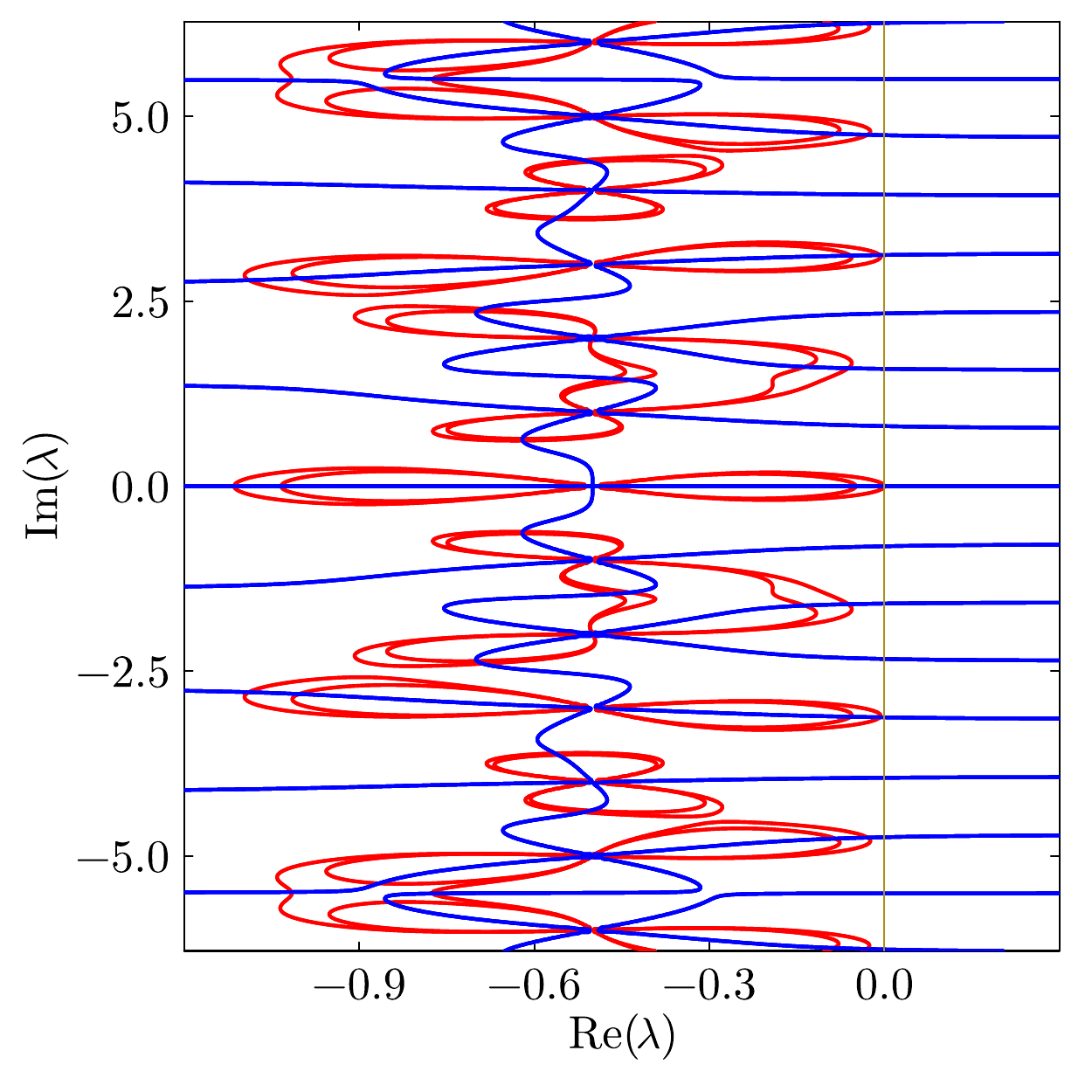} \hspace*{1cm}
\includegraphics[width=0.4\textwidth]{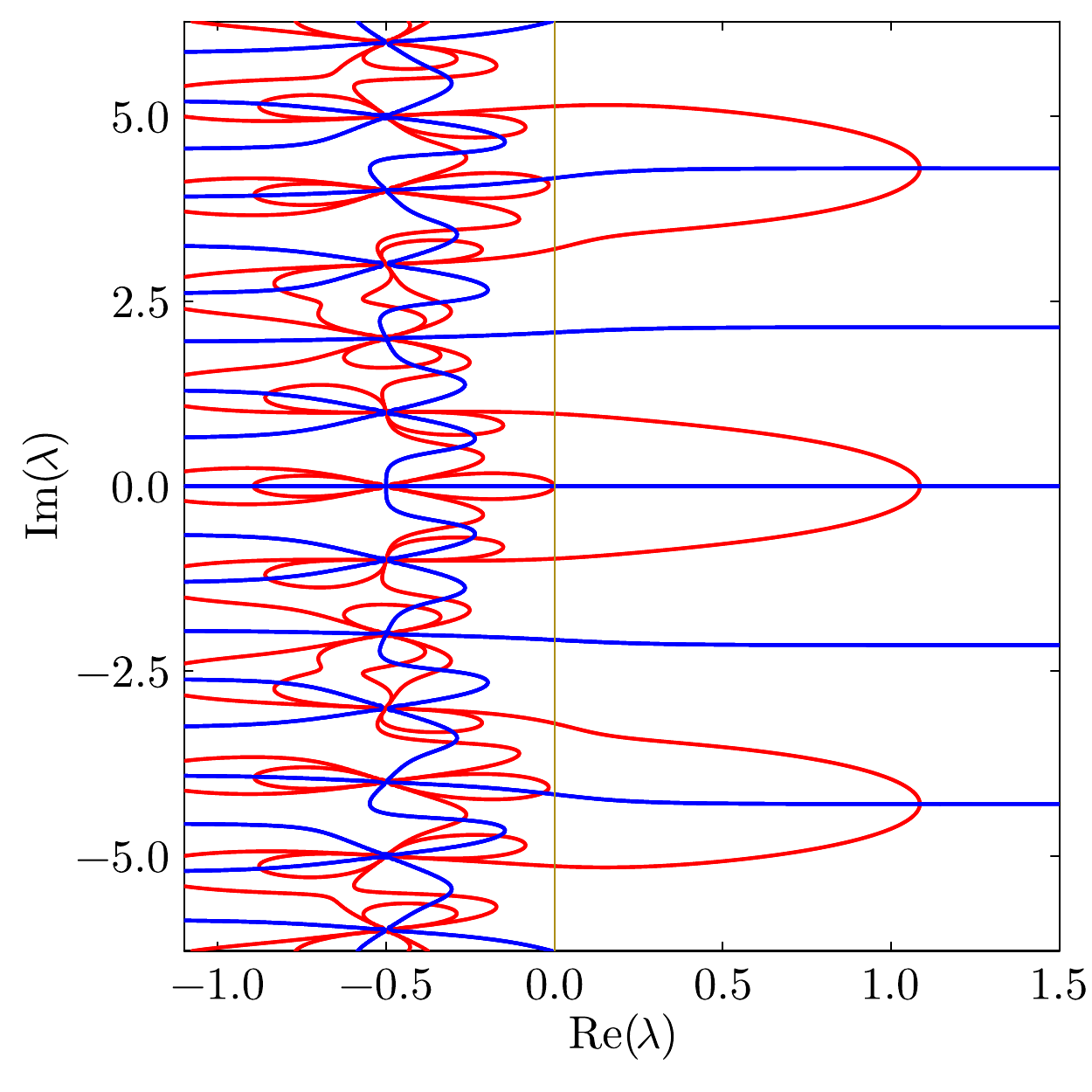} 
\caption{
\label{Fig:BumpSpectrum}
A calculation of the spectrum for the branch of solutions shown in Fig.~\ref{Fig:Deltavsrho} at $\rho=5$ and $\alpha=0.5$.
Here the red and blue curves show the zero levels sets of $\text{Re} \, \mathcal{E}_\pm(\lambda)$ and $\text{Im} \, \mathcal{E}_\pm(\lambda)$ respectively.  An eigenvalue is found wherever the red and blue lines cross.  
Left:  Spectrum on the upper branch of Fig.~\ref{Fig:Deltavsrho} left panel.  Here all the eigenvalues reside in the left had complex plane and the bump is stable.  Right:  Spectrum on the lower branch of Fig.~\ref{Fig:Deltavsrho} left panel.  Here, some eigenvalues can be found in the right hand complex plane showing that the bump is unstable.}
\end{figure}
Some example spectra are plotted in Fig.~\ref{Fig:BumpSpectrum}.  Here we plot the zero levels sets of $\text{Re} \, \mathcal{E}_\pm(\lambda)$ and $\text{Im} \, \mathcal{E}_\pm(\lambda)$ respectively.  Eigenvalues (zeros of $\mathcal{E}_\pm(\lambda)$) are found where the two sets of lines intersect.  Referring to the right hand panel of Fig.~\ref{Fig:Deltavsrho} for $\Delta = \Delta(\alpha)$ (for fixed $\rho$) we find that the lower branch of solutions is unstable.  On the upper branch the spectrum in the left hand panel of Fig.~\ref{Fig:BumpSpectrum} for $\alpha=0.5$ does not show any eigenvalues in the right hand complex plane.  However, there are eigenvalues (other than the zero from translation invariance) that are close to the imaginary axis and that more closely approach it with an increase in $\alpha$.  Even before reaching the limit point of the curve $\Delta = \Delta(\alpha)$ in Fig.~\ref{Fig:Deltavsrho} (right panel) at $\alpha \simeq 2$ it appears that the upper branch can become marginally unstable (as suggested by the location of the green circle on the upper branch).
Thus, with increasing $\alpha$ a bump solution (with a fixed pair $(\Delta, \rho)$) may go unstable or even cease to exist.  In this case we can resort to direct numerical simulations to determine emergent behaviour.

\begin{figure}
\centering
\includegraphics[width=0.4\textwidth]{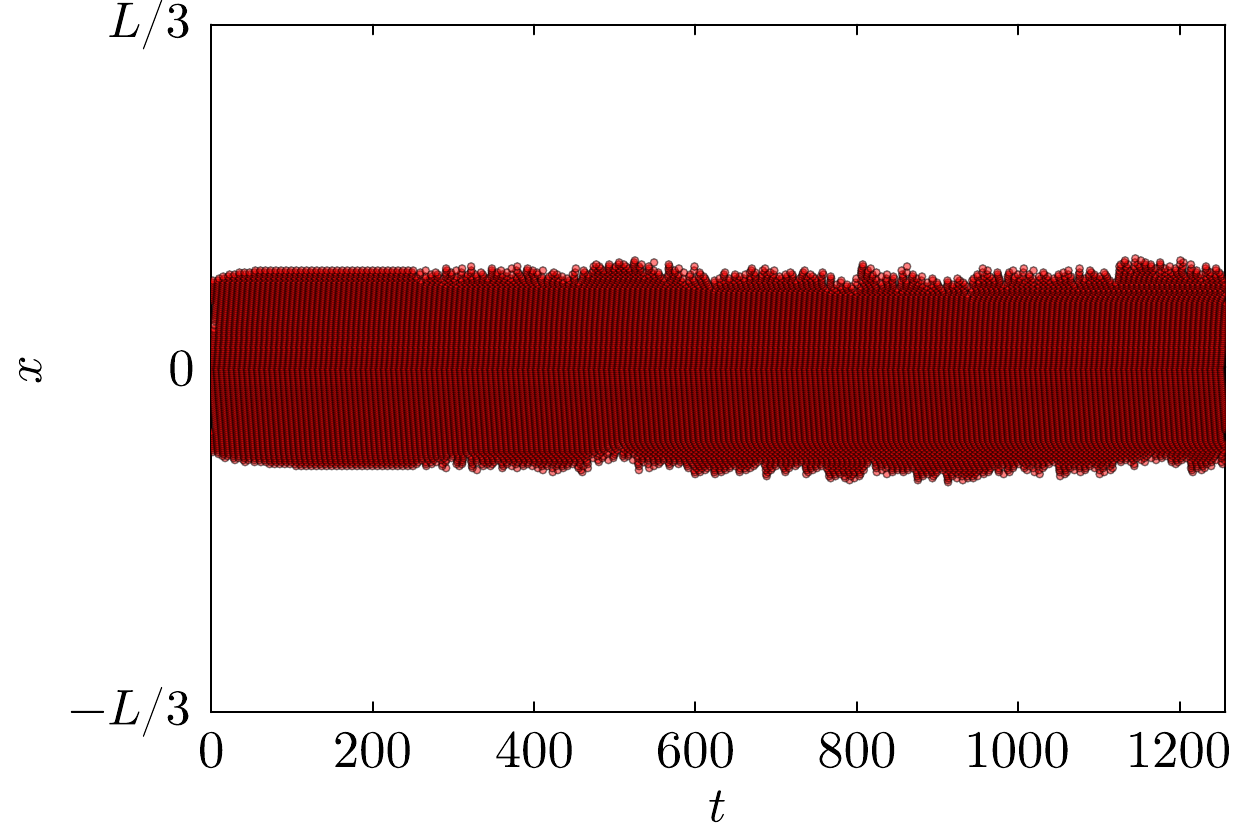} \hspace*{1cm}
\includegraphics[width=0.4\textwidth]{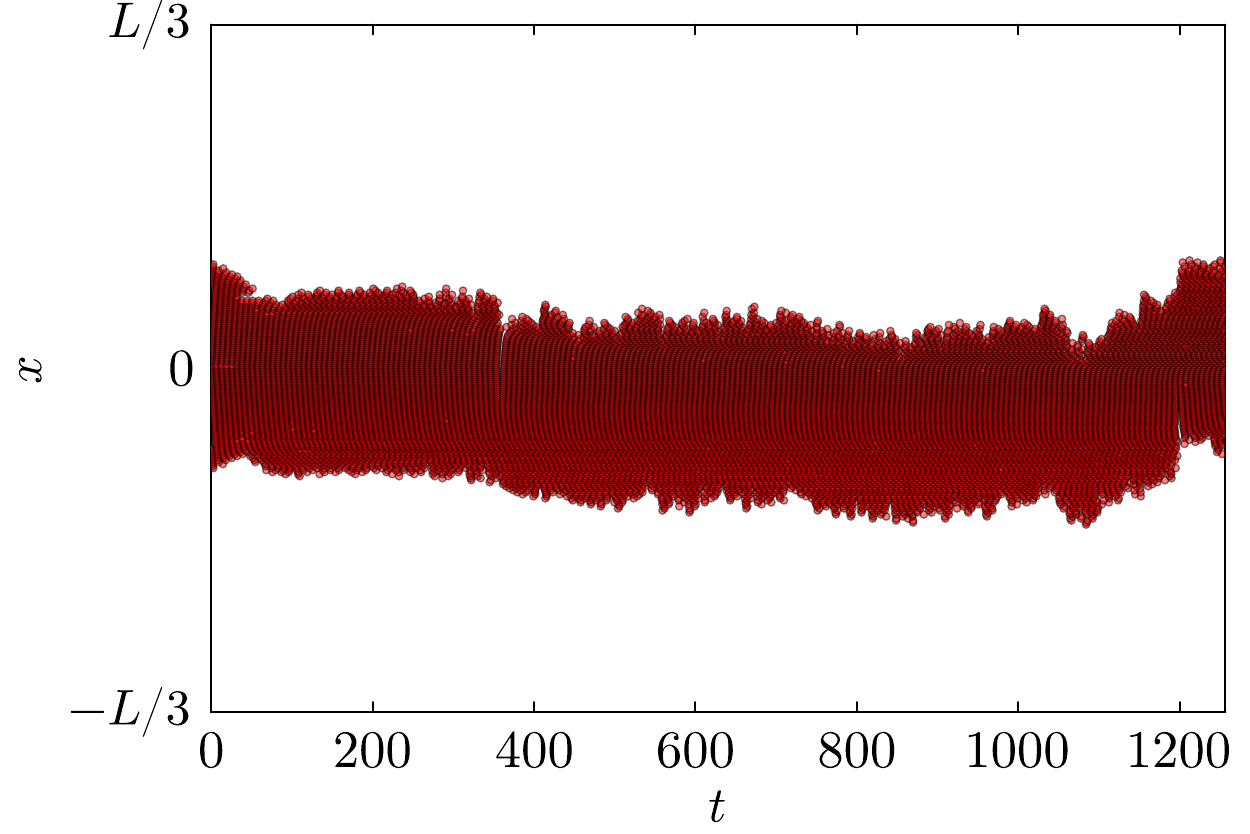}
\caption{
\label{Fig:Wandering}
Raster of a wandering bump from direct numerical simulations with parameters as in Fig.~\ref{Fig:Deltavsrho} and initiated from a bump with $\rho=5$.
Here, $\alpha$ is initially set to $\alpha = 0.5$ and one fifth of the way through the simulation it is switched to $\alpha=1$.
Simulations were performed using a regular spatial discretisation on an interval $L = 8 \sigma$ using $2^{10}$ mesh points and periodic boundary conditions.
Left:  Heaviside firing rate.  Right: Smooth firing rate with $r = 10^{-5}$.  All other parameters as in Fig.~\ref{Fig:Deltavsrho}.  
}
\end{figure}

In Fig.~\ref{Fig:Wandering} we show direct simulations of a continuum Lighthouse model that exhibit so-called \text{wandering} states.  For the choice of a Heaviside nonlinearity (left panel) we prepare the system such that a bump of the form (\ref{bump}) can persist robustly for some choice of $\alpha$.  At some point in time we then switch the value of $\alpha$ to a larger value, at which point we start to see the bump edges fluctuate and a less regular firing patter emerge.  Overall, a bump still persists though it is no longer stationary (in the sense that the same set of neurons does not continue to fire indefinitely).  Moreover, the travelling waves that make up the core of the bump are less regular than the periodic ones considered in our analysis.
Similar wandering behaviour is observed with the smooth function given by (\ref{S}) with $r$ small (right panel).  Here, we see that fluctuations in the period of neurons defining the bump edge can have a larger impact on the wandering pattern.  Recent work by Avitabile \textit{et al}. for IF networks has suggested that 
wandering bump attractors manifest a form of turbulence, and that the building blocks of the bump attractor can be understood in terms of (echoes of) unstable waves \cite{Avitabile2023}.  The same notions apply equally well to the wandering bumps of the Lighthouse model.

\section{Discussion\label{Sec:Discussion}}

In this paper we have revisited the Lighthouse model of Hermann Haken.  Building on his legacy we have presented new tools for the analysis of waves, bumps, and patterns that highlight the usefulness of this spiking model in theoretical neuroscience. Indeed, it seems ideally suited for further study to gain insight into patterns seen in balanced cortical networks with clustered connections \cite{Kumar2012}, and for theorising about \textit{metastability} in spiking systems (patterns of firing activity across simultaneously recorded neurons that linger for 300 ms to 3 s prior to transitioning to a new pattern) \cite{Brinkman2022}.  To pursue this it would be very advantageous to relax the convenient row sum constraint that we have adopted here in our study of discrete networks, and the assumption of translation invariance in the continuum model, and treat more heterogeneous systems and solutions \cite{Spreizer2019}.  This heterogeneity could also be included at the node level with the introduction of threshold noise \cite{Thul2016}.  Furthermore, although we introduced multiple communication delays into the network model, we only showed how to treat a single common delay.  Finally, it is well to remember that we have focused mainly on spontaneous activity and that developing results for driven spiking networks is very important for understanding cortical dynamics \cite{Doiron2014} and the observation that natural and pharmacologically driven brain state transitions primarily impact the spontaneous rate of slow-firing neurons \cite{Dearnley2023}.  

As well as providing insight into biological networks the computational simplicity of the Lighthouse model lends itself to uses in the growing field of bio-inspired AI and neuromorphic engineering \cite{Guo2023,DiCaterina2024}, particularly in areas where energy efficiency and biologically plausible computation are important.  Altogether, the Lighthouse model offers an exciting and fertile ground for discovery, bridging from theoretical neuroscience to new potential applications in AI that make use of spiking neural networks.

\section*{Appendix A\label{Appendix:Saltation}}

For a recent discussion of saltation operators for nonsmooth dynamical systems see \cite{Coombes2023}[Ch.~ 3].
In the context of this paper, the saltation operator of Sec.~\ref{sec:MSF} is a $3 \times 3$ matrix in the form
\begin{equation} K(T) = \D g(x(T^-)) +   \frac{\left [\dot{x}(T^+)
  - \D g(x(T^-)) \dot{x}(T^-) \right ] \left [ \nabla_{x}
  h(x(T^-)) \right ]^\mathsf{T}}{\nabla_{x} h(x(T^-)) \cdot \dot{x}
(T^-)} .  
\label{Saltation} 
\end{equation}
For the linearised Lighthouse model we have that $\D g = I_3$ (the $3 \times 3$ identity matrix), and $\nabla_{x}h = (1,0,0)$.  At a firing event the $\alpha$-function has a jump such that $(\dot{s}, \dot{u}) \rightarrow (\dot{s} +\alpha^2, \dot{u}-\alpha^2)$.

\section*{Appendix B\label{Appendix:BrainWave}}

For a recent discussion of the derivation of the brain wave equation from a non-local and delayed source term in a neural field, see \cite[Ch. 9]{Coombes2023}.
In one spatial dimension the integral $\psi(x,t) = \int_{\RSet} \d y \, E(y;\sigma) s(x-y, t-|y|/v)$, where $E(x;\sigma) = \e^{-|x|/\sigma}/(2 \sigma)$, can be constructed as the solution to the PDE:
\begin{equation}
    \left[\left(\frac{1}{\sigma}+\frac{1}{v}\PD{}{t}\right)^2-\PD{^2}{x^2}\right]\psi=\frac{1}{\sigma}\left(\frac{1}{\sigma}+\frac{1}{v}\PD{}{t}\right) s.
\end{equation}
In two spatial dimensions, the choice $E(r;\sigma) = \e^{-r/\sigma}/(2 \pi \sigma^2)$ yields the \textit{approximate} PDE
\begin{equation}
    \left[\left(\frac{1}{\sigma}+\frac{1}{v}\PD{}{t}\right)^2- \frac{3}{2} \nabla^2 \right]\psi= \sigma^{-2}s.
\end{equation}
This PDE is approximate in the sense that it is only expected to accurately capture the behaviour of long wavelength patterns.  However, it can be viewed as an exact description for the choice\\
$E(r;\sigma) =  K_0(\sqrt{2/3} r|/\sigma)/(3 \pi \sigma^2)$, where $K_0$ is the modified Bessel function of the second kind of order zero.

\section*{Appendix C\label{Appendix:Efunction}}

Introduce the integral representation
\begin{equation}
\eta'(t) = \frac{1}{2 \pi} \int_{-\infty}^ \infty \d z \, (\i z)\widehat{\eta}(z) \e^{\i z t}.
\end{equation}
For $u_m(x) = \e^{\i k x} \e^{\lambda [m + \rho x/T]}$ equation (\ref{A}) takes the form
\begin{align}
\sum_{p \in \ZSet}
\int_0^T \d t \, S'(\psi(t;T,\rho)) \int_{\RSet} \d y \, w(y)  
\int_{-\infty}^\infty \frac{\d z}{2 \pi}  (\i z)\widehat{\eta}(z) \e^{\i z (t + pT + \rho y- |y|/v)]}
\e^{\i k (x-y)} \e^{\lambda [m-p + \rho (x-y)/T]}.
\end{align}
The summation over $p$ can now be performed using the Dirac comb result that 
$\sum_p \e^{\i p (zT +\i \lambda)} = 2 \pi/T \sum_n \delta(z - \i \lambda /T +\omega_n)$, where $\omega_n = 2 \pi n /T$.  Collecting terms gives a reduction to the form $a(\lambda,k;\rho) u_m(x)$, where
\begin{align}
a(\lambda, k;\rho) &=  \frac{1}{T} \sum_{n \in \ZSet}
\int_0^T \d t \, S'(\psi(t;T,\rho)) 
\e^{\i (\omega_n -\i \lambda/T)} (\i) (\omega_n - \i \lambda/T) \times \nonumber \\
& \widehat{\eta} (\omega_n - \i \lambda/T)\int_{\RSet} \d y \, w(y)  
\e^{- \i y (k-\omega_n \rho)} \e^{-\i |y| (\omega_n - \i \lambda/T)/v} ,
\end{align}
from which we obtain (\ref{arho}).

\bibliographystyle{unsrt}
\bibliography{HakenlightHouse}

\end{document}